\documentstyle[a4,amstex,amscd,euscript,epsbox,amssymb]{amsart}
\numberwithin{equation}{section}
\newtheorem{thm}{Theorem}[section]
\newtheorem{cor}[thm]{Corollary}
\newtheorem{prop}[thm]{Proposition}
\newtheorem{lem}[thm]{Lemma}
\theoremstyle{remark}
\newtheorem{rem}{Remark}[section]
\newtheorem{exam}{Example}[section]
\begin{document}
\newcommand{\dcp}{\,{\triangleright \kern-0.145em \triangleleft}}
\newcommand{\dcpP}{\,{\triangleright \kern-0.145em \triangleleft}_{\cal P}}
\newcommand{\tensE}{{\otimes}_{\frak E}}
\newcommand{\tenc}{{\intercal}}
\newcommand{\cu}{\varepsilon}
\newcommand{\enal}{e_{\lambda}}
\newcommand{\enam}{e_{\mu}}
\newcommand{\enan}{e_{\nu}}
\newcommand{\enax}{e_{\xi}}
\newcommand{\ema}{{\stackrel{\scriptscriptstyle\circ}{e}}}
\newcommand{\emal}{{\stackrel{\scriptscriptstyle\circ}{e}_{\lambda}}}
\newcommand{\emam}{{\stackrel{\scriptscriptstyle\circ}{e}_{\mu}}}
\newcommand{\eman}{{\stackrel{\scriptscriptstyle\circ}{e}_{\nu}}}
\newcommand{\emax}{{\stackrel{\scriptscriptstyle\circ}{e}_{\xi}}}
\newcommand{\emai}{{\stackrel{\scriptscriptstyle\circ}{e}_i}}
\newcommand{\wsta}{{\stackrel{\scriptscriptstyle *}{w}}}
\newcommand{\ehoi}{e_{{\HH}^{\circ},i}}
\newcommand{\ehoj}{e_{{\HH}^{\circ},j}}
\newcommand{\ehok}{e_{{\HH}^{\circ},k}}
\newcommand{\ehol}{e_{{\HH}^{\circ},l}}
\newcommand{\emahoi}{{\stackrel{\scriptscriptstyle\circ}{e}_{{\HH}^{\circ},i}}}
\newcommand{\emahoj}{{\stackrel{\scriptscriptstyle\circ}{e}_{{\HH}^{\circ},j}}}
\newcommand{\emahok}{{\stackrel{\scriptscriptstyle\circ}{e}_{{\HH}^{\circ},k}}}
\newcommand{\emahol}{{\stackrel{\scriptscriptstyle\circ}{e}_{{\HH}^{\circ},l}}}
\newcommand{\epma}{{\stackrel{\scriptscriptstyle\circ}{\varepsilon}}}
\newcommand{\epmai}{{\stackrel{\scriptscriptstyle\circ}{\varepsilon}_i}}
\newcommand{\epmaj}{{\stackrel{\scriptscriptstyle\circ}{\varepsilon}_j}}
\newcommand{\epmak}{{\stackrel{\scriptscriptstyle\circ}{\varepsilon}_k}}
\newcommand{\epmal}{{\stackrel{\scriptscriptstyle\circ}{\varepsilon}_l}}
\newcommand{\Hom}{\mathrm{Hom}}
\newcommand{\End}{\mathrm{End}}
\newcommand{\Tr}{\mathrm{Tr}}
\newcommand{\id}{\mathrm{id}}
\newcommand{\card}{\mathrm{card}}
\newcommand{\spa}{\mathrm{span}}
\newcommand{\op}{\mathrm{op}}
\newcommand{\cop}{\mathrm{cop}}
\newcommand{\bop}{\mathrm{bop}}
\newcommand{\HH}{\frak H}
\newcommand{\KK}{\frak K}
\newcommand{\Hhat}{\hat{\frak H}}
\newcommand{\Ht}{\tilde{\frak H}}
\newcommand{\HGm}{\frak{H} [G^{-1}]}
\newcommand{\Hc}{\mathrm{Hc}}
\newcommand{\HcH}{\mathrm{Hc}({\frak H})}
\newcommand{\Aw}{{\frak A}(w)}
\newcommand{\SS}{{\frak S}(A_{N-1};t)_{\epsilon}}
\newcommand{\SSe}{{\frak S}(A_{N-1};t)_{\epsilon,\zeta}}
\newcommand{\SSz}{{\frak S}(A_{N-1};t)_{\epsilon,\zeta}}
\newcommand{\SSi}{{\frak S}(A_{N-1};t)_{\epsilon,\zeta}^{\iota}}
\newcommand{\Ss}{{\frak S}(A_{1};t)_{\epsilon}}
\newcommand{\AS}{{\frak A}(w_{N,t,\epsilon})}
\newcommand{\ASe}{{\frak A}(w_{N,t,\epsilon,\zeta})}
\newcommand{\Str}[1]{\mathrm{Str}^{#1} (\G, *)}
\newcommand{\Len}{\EuScript{L}}
\newcommand{\sgn}[1]{(- \epsilon)^{\EuScript{L} ( #1 )}}
\newcommand{\Rp}{{\cal R}^+}
\newcommand{\Rm}{{\cal R}^-}
\newcommand{\Rpm}{{\cal R}^{\pm}}
\newcommand{\Rmp}{{\cal R}^{\mp}}
\newcommand{\Rph}{{\hat{\cal R}}^+}
\newcommand{\Rmh}{{\hat{\cal R}}^-}
\newcommand{\Rpmh}{{\hat{\cal R}}^{\pm}}
\newcommand{\Rmph}{{\hat{\cal R}}^{\mp}}
\newcommand{\Rt}{{\tilde{\cal R}}}
\newcommand{\Rpt}{{\tilde{\cal R}}^+}
\newcommand{\Rmt}{{\tilde{\cal R}}^-}
\newcommand{\Rpmt}{{\tilde{\cal R}}^{\pm}}
\newcommand{\G}{\EuScript G}
\newcommand{\V}{\EuScript V}
\newcommand{\K}{\Bbb K}
\newcommand{\C}{\Bbb C}
\newcommand{\R}{\Bbb R}
\newcommand{\Q}{\Bbb Q}
\newcommand{\Z}{\Bbb Z}
\newcommand{\LD}{\mathrm{LD}}
\newcommand{\GLD}{\G_{\mathrm{LD}}}
\newcommand{\wLD}{w_{\mathrm{LD}}}
\newcommand{\KP}{\mathrm{KP}}
\newcommand{\GLEH}{\mathrm{GLE} (\HH)}
\newcommand{\LLam}{\boldsymbol\Lambda}
\newcommand{\aaa}{\bold a}
\newcommand{\bbb}{\bold b}
\newcommand{\ccc}{\bold c}
\newcommand{\ddd}{\bold d}
\newcommand{\p}{\bold p}
\newcommand{\q}{\bold q}
\newcommand{\r}{\bold r}
\newcommand{\s}{\bold s}
\newcommand{\taaa}{\tilde{\bold a}}
\newcommand{\tbbb}{\tilde{\bold b}}
\newcommand{\tccc}{\tilde{\bold c}}
\newcommand{\tddd}{\tilde{\bold d}}
\newcommand{\tp}{\tilde{\bold p}}
\newcommand{\tq}{\tilde{\bold q}}
\newcommand{\tr}{\tilde{\bold r}}
\newcommand{\ts}{\tilde{\bold s}}
\newcommand{\st}{\frak {s}}
\newcommand{\en}{\frak {r}}
\newcommand{\suma}{\sum_{(a)}}
\newcommand{\sumb}{\sum_{(b)}}
\newcommand{\sumc}{\sum_{(c)}}
\newcommand{\sumd}{\sum_{(d)}}
\newcommand{\sumx}{\sum_{(x)}}
\newcommand{\sumy}{\sum_{(y)}}
\newcommand{\sumab}{\sum_{(a),(b)}}
\newcommand{\sumbc}{\sum_{(b),(c)}}
\newcommand{\sumad}{\sum_{(a),(d)}}
\newcommand{\sumxy}{\sum_{(x),(y)}}
\newcommand{\sumax}{\sum_{(a),(x)}}
\newcommand{\sumbx}{\sum_{(b),(x)}}
\newcommand{\sumabc}{\sum_{(a),(b),(c)}}
\newcommand{\sumabcd}{\sum_{(a),(b),(c),(d)}}
\newcommand{\sumk}{\sum_{k \in \V}}
\newcommand{\suml}{\sum_{l \in \V}}
\newcommand{\sumkl}{\sum_{k,l \in \V}}
\newcommand{\face}[4]{\left( \scriptstyle{#1} 
                      \textstyle{\frac[0pt]{#2}{#3}} \scriptstyle{#4} \right)}
\newcommand{\Wrpqs}[4]{w \! \left[ #1 \, \frac[0pt]{#2}{#3} \, #4 \right]}
\newcommand{\wrpqs}[4]{w \! \left[ \scriptstyle{#1} 
                      \textstyle{\frac[0pt]{#2}{#3}} \scriptstyle{#4} \right]}
\newcommand{\hijk}[4]{ \!\! \left[ {#1 \,\, #2} \atop {#3 \,\, #4} \right]}

\title{Face algebras and unitarity of $\text{SU(N)}_{\text{L}}$-TQFT}
\author{Takahiro Hayashi \quad}
\date{$\qquad$\\
Graduate School of Mathematics, 
Nagoya University, \\ $\quad\;$
Furo-cho, Chikusa-ku, Nagoya 464, Japan}
\maketitle

\begin{abstract}
 Using face algebras (i.e. algebras of L-operators of 
 IRF models), we construct modular tensor categories
 with positive definite inner product,  
 whose fusion rules and S-matrices are the same as
 (or slightly different from)
 those obtained by $U_q (\frak{sl}_{N})$ at roots 
 of unity.				 
 Also we obtain state-sums of 
 ABF models on framed links which give quantum 
 $SU(2)$-invariants of corresponding 3-manifolds.
\end{abstract}

\section{introduction}

As is well known, quantum groups have their origin 
in the theory of quantum inverse scattering method.
More specifically, they first appeared 
as so-called algebras of L-operators
of lattice models (of vertex type).  	
For example, the simplest quantum group $U_q(\frak{sl}(2))$
can be viewed as the algebra of L-operators of 
6-vertex model without spectral parameter.  
It seems that it is worth trying to study 
algebras of L-operators 
independently from the framework of Drinfeld-Jimbo algebra.

By investigating algebraic structure of 
lattice models of face type, 
we found a new class of quantum groups, which is called
the class of {\em face algebras}
(cf. \cite{subf}-\cite{cpcq}
and also \cite{BS,JurcoSchupp,Schauenburgface}).
It contains all bialgebras as a subclass.
Moreover, as well as bialgebras, face algebras
produce monoidal categories as their (co-)module categories.

In this paper, we give a detailed study of face algebras 
$\SS$, which are obtained as algebras of L-operators of 
RSOS models of type $A_{N-1}$ $(N \geq 2)$ (cf. \cite{JMO}),		
where $\epsilon = \pm 1$ and $t$ denotes
a primitive $2(N + L)$-th root of unity with $L \geq 1$.
We also give two applications of $\SS$ to 3-dimensional
topological quantum field theory (TQFT) and 
corresponding quantum invariants of 3-manifolds.  

We show that the algebra $\SS$ is finite-dimensional cosemisimple
and that its dual is a $C^*$-algebra for a suitable $t$. 
Also, we classify irreducible comodules of $\SS$ and determine
their dimensions.
Moreover, we construct various structures on $\SS$,
such as the antipode, the braiding and the ribbon functional.  
The algebra $\SS$ is constructed as a quotient
of the face version $\AS$ of FRT construction modulo 
one additional relation ``$\det = 1$'',
where $w_{N,t,\epsilon}$ is the Boltzmann weight of RSOS models  
of type $A_{N-1}$ without the spectral parameter 
and $\det$ denotes the ``(quantum) determinant'' of 
$\AS$. 
Since the representation theory of $\AS$ is relatively 
easily established using a result on Iwahori-Hecke 
algebras due to H. Wenzl, the core of our work is to 
study the properties of the element 
$\det$ or, corresponding ``exterior'' algebra. 

Next, we explain the unitarity of 3-dimensional TQFT briefly.
Roughly speaking, a 3-dimensional TQFT is a map which
assigns to each 3-cobordism $(M, \partial_- M,$ $\partial_+ M)$, 
a linear map 				
$\tau (M)\!: \EuScript{T}(\partial_- M) \to \EuScript{T}(\partial_+ M)$.
Here, by a 3-dimensional cobordism, we mean 
a compact 3-dimensional manifold $M$ 
whose boundary is a disjoint union of
two closed surfaces $\partial_- M, \partial_+ M$. 
A 3-dimensional TQFT is called unitary if 
$\EuScript{T}(\partial_{\pm} M)$ are (finite-dimensional)
Hilbert spaces and $\tau (-M) = \tau (M)^*$ 
for each 3-cobordism $(M, \partial_- M, \partial_+ M)$.
It is established in \cite{Turaev} that to obtain 
a (unitary) 3-dimensional MTC, 
it suffices to construct a (unitary) modular tensor category 
(MTC) (i.e. a braided category which satisfies
certain additional properties).				
The most important examples of MTC are constructed as 
semisimple quotients $\cal{C} (\frak{g}, q)$ 
of some module categories of Drinfeld-Jimbo algebra 
$U_q (\frak{g})$ at roots $q$ of unity
(cf. \cite{ReshetikhinTuraev,Andersen,GelfandKazhdan,Kirillov},
and see also \cite{Finkelberg,TuraevWenzl,TuraevWenzl2}
for other construction of MTC's).
For a suitable $q$, it also is expected that 
$\cal{C} (\frak{g}, q)$ is a unitary MTC.
However, it seems that it is not easy to verify it directly, 
since $U_q (\frak{g})$ is non-semisimple and cannot have a 
$C^*$-algebra structure (cf. A. Kirillov, Jr, \cite{Kirillov} and 
V. Tuaev and H. Wenzl \cite{TuraevWenzl2}).

The first application of $\SS$ is to show that 
the category
$\cal{C}_{\frak{S}}(A_{N-1},t)_{\epsilon}$ of 
all finite-dimensional right $\SS$-comodules is a MTC,
and that (for a suitable $t$), the category
$\cal{C}^u_{\frak{S}}(A_{N-1},t)_{\epsilon}$ of 
all finite-dimensional unitary
$\SS$-comodules is a unitary MTC,
whose fusion rules agree with 
those of $\cal{C} (\frak{sl}_N, q)$.
Here ``unitary'' comodule means a comodule with inner product
which satisfies some conditions.
The quantum dimensions and $S$-matrices of
$\cal{C}_{\frak{S}}(A_{N-1},t)_{\epsilon}$ are the same as
(or slightly different from) those of
$\cal{C} (\frak{sl}_N, q)$  
(according to the choice of $\epsilon$ and 
another sign parameter $\iota$ when $N$ is even).
Although we use $U_q (\frak{sl}_N)$ to obtain some
combinatorial formulas, the essential part of our theory
is independent from $U_q (\frak{sl}_N)$.
Hence the equivalence of $\cal{C} (\frak{sl}_N, q)$ and
$\cal{C}_{\frak{S}}(A_{N-1},t)_{\epsilon}$ 
is left as an open problem.
However, by the result of Kazhdan and Wenzl \cite{KazhdanWenzl},
these two categories are equivalent up to a ``twist.''

Unlike the module category of $U_q(\frak{g})$, 
the category  
$\cal{C}_{\frak{S}}(A_{N-1},t)_{\epsilon}$ 
itself is semisimple.
Moreover, it has a apparent similarity to
the spaces of the conformal blocks of 
Wess-Zumino-Witten (WZW) models.
We hope that there exists a direct connection
between $\SS$ and WZW models,
similarly to Drinfeld \cite{DrinfeldGal}.

The second application is to give an explicit description of 
the quantum $SU(2)$-invariant $\tau (M)$ 
of closed 3-manifolds $M$,
which is associated with
$\cal{C}_{\frak{S}}(A_1,t)_{\epsilon}$.
More precisely, we express $\tau (M)$
as a state sum on each generic link diagram $D$
which represents $M$.
It gives a direct connection between the invariant
and ABF model.

The paper is organized as follows.
We start in  Sect. 2-4, by recalling elementary properties of 
face algebras $\HH$, various structures on $\HH$
and their relations to the comodule category 
of $\HH$.
In Sect. 5, we recall the notion of star-triangular (Yang-Baxter)
face models and flat face models. The later is introduced in
\cite{gal}, and is a variant of Ocneanu's notion of 
flat biunitary connection in operator algebra.
Flat face models play an crucial role to the determination 
of the representation theory of $\SS$. 
In Sect. 6-7, we define the algebras $\SS$ and state the main result
on the representation theory of them.
In Sect. 9-10, we construct several structures on
$\SS$, which we call the transpose, the costar structure, 
the antipode and the ribbon functional.
Consequently, we see that 
$\cal{C}_{\frak{S}}(A_{N-1},t)_{\epsilon}$
(resp. $\cal{C}^u_{\frak{S}}(A_{N-1},t)_{\epsilon}$)
is a (unitary) ribbon category.
In Sect. 11, we prove that 
$\cal{C}_{\frak{S}}(A_{N-1},t)_{\epsilon}$
is a modular tensor category, by computing its $S$-matrix.
Section 8 and Section 12 are devoted to some technical 
calculations on face analogues of the exterior algebras. 
In Sect. 13-14, we give an explicit description of 
$\cal{C}_{\frak{S}}(A_1,t)_{\epsilon}$ 
and give a state sum expression of 
the quantum $SU(2)$-invariant stated above. \par\noindent 

After submitting the manuscript, one of the referee informed me
the existence of the following two papers;\par\noindent 
H. Wenzl, $C^*$-tensor categories from quantum groups, 
{\it J. Amer. Math. Soc.} {\bf 11} (1988), pp. 261-282.\par\noindent 
C. Blanchet, Heck algebras, modular categories and 3-manifolds 
quantum invariants, preprint.
				
The former gives a proof of the unitarity of 
$\cal{C} (\frak{g}, q)$ for each $\frak{g}$.
The latter gives a construction of modular tensor categories via 
Iwahori-Hecke algebras at root of unity.

Throughout this paper, we use Sweedler's sigma notation for 
coalgebras, such as 					
$(\Delta \otimes \mathrm{id}) (\Delta (a)) =
\sum_{(a)} a_{(1)} \otimes a_{(2)} \otimes a_{(3)}$
(cf. \cite{Sweedler}).

\section{face algebras}

%
In this section, we give the definition of 
the face algebra and various structures on it.
Let $\HH$ be an algebra over a field ${\Bbb K}$ equipped
with a coalgebra structure $(\HH,\Delta,\cu)$.  
Let ${\EuScript V}$ be a finite non-empty set and 
let $\enal$ and $\emal$
$(\lambda \in {\V})$ be elements of ${\HH}$. 
We say that $({\HH}, {\enal,\emal})$ is a {\em ${\V}$-face algebra} 
if the following relations are satisfied:
\begin{equation}
 \Delta (ab) = \Delta(a) \Delta(b),
\label{D(ab)}
\end{equation}
\begin{equation}
 \enal \enam = \delta_{\lambda \mu} \enal ,
 \quad \emal \emam = \delta_{\lambda \mu} \emal,
 \quad \enal \emam = \emam \enal ,
\label{ee}
\end{equation}
\begin{equation}
 \sum_{\nu \in {\V}} \, \enan = 1 = \sum_{\nu \in \V} {\eman},
\label{sume}					
\end{equation}
\begin{equation}			
 \Delta(\emal \enam) = \sum_{\nu \in \V} \emal \enan \otimes \eman \enam,
\quad
 \cu (\emal \enam) = \delta_{\lambda \mu} ,
\label{D(ee)}
\end{equation}
\begin{equation}
 \cu (ab) = \sum_{\nu \in \V} \cu (a\enan) \cu (\eman b)
\label{e(ab)}
\end{equation}
for each $a, b \in \HH$ and $\lambda, \mu \in \V$.
We call elements $\enal $ and $\emal$ {\it face idempotents}
of $\HH$. 
It is known that bialgebra is an equivalent notion of 
$\V$-face algebra with $\mathrm{card}(\V) = 1$.
For a $\V$-face algebra $\HH$, we have the following formulas:
\begin{equation}
 \cu (\emal a) = \cu (\enal  a), 
 \quad \cu (a \emal) = \cu (a \enal ), 
\label{e(eae)}
\end{equation}
\begin{equation}
 \suma a_{(1)} \cu  (\enal  a_{(2)} \enam) = \enal  a \enam,  
\label{ae(eae)}
\end{equation} 
\begin{equation}
 \suma \cu (\enal  a_{(1)} \enam)a_{(2)} = \emal a \emam,  
\label{e(eae)a}
\end{equation}
\begin{equation}
 \Delta (a) = \sum_{\nu, \xi} \suma 
\enan a_{(1)} \enax \otimes \eman a_{(2)} \emax,  
\label{D(a)}
\end{equation}
\begin{equation}
 \suma \enal  a_{(1)} \enam \otimes a_{(2)} = 
 \suma a_{(1)} \otimes \emal a_{(2)} \emam,  
\label{eae*a}
\end{equation}
\begin{equation} 			
 \Delta (\emal \enam a \ema_{\lambda'} e_{\mu'}) =
 \suma \emal a_{(1)} \ema_{\lambda'} \otimes \enam a_{(2)} e_{\mu'}  
\label{D(eeaee)}
\end{equation}
for each $a \in \HH$ and $\lambda, \mu, \lambda', \mu' \in \V$.


%
%
%
%
%
%
Let $\G$ be a finite oriented graph 
with set of vertices ${\V}$ = ${{\G}^0}$.
For an edge ${\p}$, we denote by ${\st (\p)}$ and ${\en (\p)}$
its {\em source} ({\em start}) and its {\em range} ({\em end}) respectively.
For each $m \geq 1$, we denote by 
${{\G}^m = {\coprod}_{\lambda, \mu \in \V}{\G}_{\lambda \mu}^m}$ 
the set of {\em paths} of $\G$ of {\em length} $m$, that is, 
$\p \in {\G}_{\lambda \mu}^m$ if $\p$ is a sequence
$(\p_1, \ldots, \p_m)$ of edges of $\G$ such that
$\st (\p):= \st (\p_1) = \lambda$,
$\en (\p_n) = \st(\p_{n+1})\,$ $(1 \leq n < m)$ 
and $\en (\p):= \en(\p_m) = \mu$.
Let ${\HH (\G)}$ be the linear span of 
the symbols ${e{\p \choose \q}}$
\( ( \p,\q \in {\G}^m, m \geq 0) \).
Then ${\HH (\G)}$ becomes a $\V$-face algebra by setting
\begin{equation}
 \emal = \sum_{\mu \in \V}e{\lambda \choose \mu},
 \quad \enam = \sum_{\lambda \in \V}e{\lambda \choose \mu},
\label{eHGi}
\end{equation}
\begin{equation}
 e{\p \choose \q} e{\r \choose \s} =
 {\delta}_{\en(\p) \st(\r)} \, {\delta}_{\en(\q) \st(\s)}
 \: e{\p\cdot\r \choose \q\cdot\s},
\label{epqers}
\end{equation}
\begin{equation}
 \Delta \left( e{\p \choose \q} \right) =
 \sum_{{\bold t} \in {\EuScript G}^m}
 e{\p \choose {\bold t}} \otimes e{{\bold t} \choose \q},
 \quad
 \cu \left( e{\p \choose \q} \right)  = {\delta}_{\p \q}
\label{D(epq)}
\end{equation}
for each $\p,\q \in {\G}^m$ 
and $ \r,\s \in {\G}^n$ $( m,n \geq 0)$. 
Here for paths $\p = (\p_1,\ldots,\p_m)$ and 
$\r = (\r_1,\ldots,\r_n)$, we set 
$\p \cdot \r = (\p_1,\ldots,\p_m,\r_1,\ldots,\r_n)$ 
if $\en(\p) = \st(\r)$.
We note that
$\K \G^m = \bigoplus_{\p \in \G^m} \K \p$
$(m \geq 0)$ becomes a right $\HH (\G)$-comodule via 
\begin{equation}
 \q \mapsto
 \sum_{\p \in \G^m}
 \p \otimes e{\p \choose \q}.
\end{equation}
\begin{prop}[\cite{fa}]
 Every finitely generated $\V$-face algebra is isomorphic to  
 a quotient of ${\HH (\G)}$ for some $\G$.
\end{prop}
We say that a linear map $S\!: {\HH} \to {\HH}$ is an {\em antipode} of ${\HH}$, 
or $({\HH}, S)$ is a {\em Hopf} ${\V}$-{\em face algebra} if 
\begin{equation}		
 \suma S(a_{(1)})a_{(2)} = \sum_{\nu \in \V} \cu (a \enan) \enan, 
 \quad				
 \suma a_{(1)}S(a_{(2)}) = \sum_{\nu \in \V} \cu (\enan a){\eman},
 \label{S(a)a}
\end{equation}
\begin{equation}
 \suma S(a_{(1)})a_{(2)}S(a_{(3)}) = S (a) 
 \label{S(a)aS(a)}
\end{equation} 
for each $a \in \HH$.
An antipode of a $\V$-face algebra is 
an antialgebra-anticoalgebra map, which satisfies
\begin{equation}
 S(\emal \enam) = \emam \enal 
 \quad (\lambda, \mu \in {\V}).
\label{S(ee)}
\end{equation}
%
%
%
The antipode of a $\V$-face algebra is unique 
if it exists.

%
%
%
%
%
%
Let ${\HH}$ be a ${\V}$-face algebra with product $m$ and 
let $\Rp$ and $\Rm$ be elements of $({\HH} \otimes{\HH})^*$.
We say that $({\HH},\Rpm)$ is a {\em coquasitriangular}
(or {\em CQT}) ${\V}$-face algebra 
if the following relations are satisfied 
\begin{equation}
 \Rp m^*(1) = \Rp, 
 \quad 
 m^*(1) \Rm = \Rm,
\label{Rpm(1)}
\end{equation}
\begin{equation}
 \Rm \Rp = m^*(1),
 \quad 
 \Rp \Rm = (m^{\text{op}})^*(1),
\label{RmRp}
\end{equation}
\begin{equation}
 \quad
 \Rp m^*(X) \Rm = (m^{\text{op}})^*(X)
 \quad (X \in \HH^*),
\label{RmXR}
\end{equation}
\begin{equation}
\label{mid(R)}
 (m {\otimes}\text{id})^* (\Rp) = \Rp_{13} \Rp_{23}, 
 \quad
 (\text{id}{\otimes}m)^* (\Rp) = \Rp_{13} \Rp_{12}.
\end{equation}
Here, as usual, for each $Z \in (\HH \otimes \HH)^*$
and $\{ i,j,k \} = \{1,2,3 \}$, we define
$Z_{ij} \in (\HH^{\otimes 3})^*$ by
$Z_{ij} (a_1, a_2, a_3) = Z(a_i, a_j) \cu (a_k)$  
$(a_1, a_2, a_3 \in \HH)$. We note, for example, that
the second formula of \eqref{mid(R)} is equivalent to
\begin{equation}
 \label{Rp(a,bc)}
 \Rp(a,\, bc) = \suma \Rp(a_{(1)},\, c) \Rp(a_{(2)},\, b)
 \quad (a, b, c \in \HH).
\end{equation}
It is known that $\Rpm$ satisfies
\begin{equation}
\label{mid(Rm)}
 (m {\otimes}\text{id})^* (\Rm) = \Rm_{23} \Rm_{13}, 
 \quad
 (\text{id}{\otimes}m)^* (\Rm) = \Rm_{12} \Rm_{13}, 
\end{equation}
\begin{equation}
 \label{Rpeeee}
 \Rp (\emal \enam a \eman \enax,\, b) = \Rp (a,\, \emam \enax b \emal \enan),
\end{equation}						
\begin{equation}
 \label{Rmeeee}						
 \Rm (\emal \enam a \eman \enax,\, b) = \Rm (a,\, \eman \enal  b \emax \enam),
\end{equation}
\begin{equation}
 \label{Rpeea}
 \Rp (\emal \enam,\, a) = \cu (\enam a \enal ),
 \quad 
 \Rp (a,\, \emal \enam) = \cu (\enal  a \enam),
\end{equation}
\begin{equation}
\label{Rmeea}
 \Rm (\emal \enam,\, a) = \cu(\enal  a \enam),
 \quad 
 \Rm (a,\, \emal \enam) = \cu(\enam a \enal )
\end{equation}
for each $\lambda, \mu, \nu, \xi \in \V$, $a,b \in \HH$.
If, in addition, $\HH$ has an antipode $S$, then
\begin{equation}
\label{Sid(R)}
 (S \otimes \text{id})^* (\Rp) = \Rm,
 \quad
 (\text{id} \otimes S)^* (\Rm) = \Rp.
\end{equation}
For a CQT Hopf ${\V}$-face algebra ${\HH}$, 
we define its {\em Drinfeld functionals} 
${\cal U}_{i} \in \HH^*$ 
$(i = 1,2)$ via
\begin{equation}
\label{Udef}
 {\cal U}_1 (a) = \suma \Rp (a_{(2)}, S(a_{(1)})),
\quad
 {\cal U}_2 (a) = \suma \Rm (S(a_{(1)}), a_{(2)})
 \quad (a \in \HH).
\end{equation}
The Drinfeld functionals 
are invertible as elements of the dual algebra
${\HH}^*$ and satisfy the following relations
\cite{gsg}:
\begin{equation}
\label{Um}
  {\cal U}_1^{-1}(a) = \suma \Rm (S(a_{(2)}),a_{(1)}),
 \quad
  {\cal U}_2^{-1}(a) = \suma \Rp (a_{(1)},S(a_{(2)})),
 \end{equation}
 \begin{equation}
 \label{U(ee)}
  {\cal U}_{i}^{\pm 1}({\emal} \enam) = {\delta}_{\lambda \mu}, 
 \quad
  {\cal U}_{i}^{\pm}(\emal a \emam) =   {\cal U}_{i}^{\pm}(\enal  a \enam), 
 \end{equation} 
 %
%
\begin{gather}
\label{m*(U)} 
 m^*({\cal U}_1) = \Rm \Rm_{21} ({\cal U}_1 \otimes {\cal U}_1)
 = ({\cal U}_1 \otimes {\cal U}_1) \Rm \Rm_{21}, \\
 m^*({\cal U}_2) 
 = \Rp_{21} \Rp ({\cal U}_2 \otimes {\cal U}_2)
 = ({\cal U}_2 \otimes {\cal U}_2 ) \Rp_{21} \Rp,
 \nonumber
\end{gather}
\begin{equation}
\label{UXU-}
 {\cal U}_i X {\cal U}_i^{-1} =(S^2)^*(X),
\end{equation}
\begin{equation}
\label{S(U)}
 S^*({\cal U}_1^{\pm 1}) = {\cal U}_2^{\mp 1}, 
\quad 
 S^*({\cal U}_2^{\pm 1}) = {\cal U}_1^{\mp 1}
\end{equation}
for each $i = 1,2$, $a \in \HH$, $X \in \HH^*$ 
and $\lambda, \mu \in \V$.
In particular, $S$ is bijective and 
${\cal U}_1{\cal U}_2^{-1}$ is a central element of ${\HH}^*$.
Let $({\HH},\Rpm)$ be a CQT Hopf $\V$-face algebra
and $\cal{V}$ an invertible central element of 
${\HH}^*$.
We say that $\cal{V}$ is a {\em ribbon functional} of $\HH$, 
or $(\HH,\cal{V})$ is a
{\em coribbon Hopf} $\V$-{\em face algebra}
if 
%
%
\begin{equation}
  m^*({\cal V}) = \Rm \Rm_{21} ({\cal V} \otimes {\cal V}),
 \label{m(V)} 
\end{equation} 
\begin{equation}
  S^* (\cal{V}) = \cal{V}.
\label{S(V)}
\end{equation}
Let $g$ (resp. ${\cal G}$) be an element of
(resp. a linear functional on) 
a $\V$-face algebra ${\HH}$. 
We say that $g$ (resp. ${\cal G}$) is {\em group-like} if 
the following relations \eqref{D(g)}-\eqref{gee}  
(resp. \eqref{G(ab)}-\eqref{G(eae)}) are satisfied:				
\begin{equation}		
\label{D(g)}					
 \Delta (g) = \sum_{\nu \in {\V}} g \enan \otimes g \eman,
\end{equation}
\begin{equation}
\label{gee}					
 g \emal \enam = \emal \enam g,
\quad
 \cu (g \emal \enam) = \delta_{\lambda \mu}, 
\end{equation}
\begin{equation}		
\label{G(ab)}					
 {\cal G}(ab) = \sum_{\nu \in {\V}}{\cal G} (a \enan ) {\cal G}({\eman} b),
\end{equation}
\begin{equation}
\label{G(eae)}					
 {\cal G}({\emal}a{\emam}) = {\cal G}(\enal a \enam),
 \quad
 {\cal G}({\emal}\enam) = \delta_{\lambda \mu}
\end{equation}
for each $a,b \in \HH$ and $\lambda, \mu \in {\V}$.
If $\HH$ has an antipode $S$, then every group-like element $g$
and group-like functional ${\cal G}$
are invertible and satisfy
\begin{equation}
 S(g) = g^{-1},
\quad
 S^*({\cal G}) = {\cal G}^{-1}. 
\end{equation} 
We denote by $\mathrm{GLE} (\HH)$ the set of 
all group-like elements of $\HH$.
\begin{prop}[\cite{rib}]
 Let ${\HH}$ be a CQT Hopf $\V$-face algebra and
 ${\cal V}$ an invertible element of $\HH^*$.  
 Then $(\HH,\cal{V})$ is a coribbon Hopf $\V$-face algebra 
 if and only if $\cal{M} = \cal{U}_1 \cal{V}^{-1}$ 
 is group-like and satisfies the following relations\rom:    
 \begin{equation}
 \label{MXM-}
 \qquad
  {\cal M} X {\cal M}^{-1} = 
  (S^2)^* (X)
  \quad (X \in \HH^*),
 \end{equation}
 \begin{equation}
 \label{M2} 
  \cal{M}^2 = \cal{U}_1 \cal{U}_2. \qquad
 \end{equation}
\end{prop}
%
%
%
For a coribbon Hopf $\V$-face algebra $(\HH, \cal{V})$, 
we call $\cal{M} = \cal{U}_1 \cal{V}^{-1}$ the 
{\it modified ribbon functional} of $\HH$
corresponding to $\cal{V}$. 
\begin{exam}
\label{L=1}
When $L =1$, the algebra $\SS$ is rather degenerate.
Hence, we treat this case here separately from the case $L \geq 2$.
(cf. \cite{wzw}).
Let $N \geq 2$ be an integer and 
$\V$ the cyclic group $\Z / N \Z$.
Let $t \in \C $ be a primitive 			
$2(N + 1)$-th root of unity and $\epsilon$ either $1$ or $-1$.
We define the $\V$-face algebra $\SS$
to be the $\C$-linear span of the symbols 
$e^i_j (m)$ $(i,j,m \in \V)$ equipped with
the structure of $\V$-face algebra given by
\begin{equation}
\label{eijpeklq} 
 e^i_j [p] e^k_l [q] = 
 \delta_{i+p, k} \delta_{j+q, l} e^i_j [p+q]
 \quad (p, q \in \Z_{\geq 0}),
\end{equation}
\begin{equation}
 \Delta( e^i_k (m) ) =
 \sum_j e^i_j (m) \otimes e^j_k (m),
\quad
 \cu (e^i_j (m)) = \delta_{ij},
\end{equation}
\begin{equation}
 \ema_i = \sum_j e^i_j (0), 
\quad
 e_j = \sum_i e^i_j (0).
\end{equation}
Here, in \eqref{eijpeklq}, we set 
\begin{equation}
 e^i_j [qN + r] = (- \epsilon)^{(i-j)q(N-1)} 
 e^i_j (r + N \Z)
\end{equation}
for each $q \in \Z$ and $0 \leq r < N$. 
The algebra $\SS$ becomes a coribbon Hopf $\V$-face algebra via  
\begin{equation}
 S (e^i_j [p]) = e^{j+p}_{i+p} [-p], 
\end{equation}
\begin{equation}
 \Rp (e^i_j [p], e^k_l [q]) = 
 \delta_{i, k+q} \delta_{jk} 
 \delta_{i+p, l+q} \delta_{j+p, l} 
 (- \zeta t)^{- pq},
\end{equation}
\begin{equation}
 \cal{V}_{\iota} (e^i_j [p]) =
 \delta_{ij} \iota^p (- \zeta t)^{p^2}, 
\end{equation}
where $i,j,k,l \in \V$, $p, q \in \Z_{\geq 0}$,
$\zeta$ denotes a solution of 
$\zeta^N = \epsilon^{N-1} t$, 
$\iota = \pm 1$ if $N \in 2 \Z$
and $\iota = 1$ if $N \in 1 + 2 \Z$.
\end{exam}

\section{comodules of face algebras}

%
In this section, we recall categorical properties of 
comodules of face algebras (cf. \cite{fa}).
We refer the readers to \cite{Kassel} for the terminologies
on monoidal (or tensor) categories.
Let $M$ be a right comodule of a $\V$-face algebra $\HH$.
We define its {\em face space decomposition} 
$M = \bigoplus_{\lambda, \mu \in {\V}}M(\lambda, \mu)$ by
\begin{equation}
 M(\lambda, \mu) = \Bigl\{ \sum_{(u)}u_{(0)} 
 \varepsilon (\enal u_{(1)} \enam) 
 \Bigm| u \in M \Bigr\}. 
\end{equation}
Here we denote the coaction $M \to M \otimes \HH$ by
$u \mapsto \sum\nolimits_{(u)}u_{(0)} \otimes u_{(1)}$
$(u \in M)$.
Let $N$ be another $\HH$-comodule. We define the 
{\em truncated tensor product} 
$M \bar{\otimes} N$ to be the $\HH$-comodule
given by
\begin{equation}
\label{MNlm}
 (M \bar{\otimes} N) (\lambda, \mu) = 
 {\bigoplus}_{\nu \in {\V}}
 M(\lambda, \nu) \otimes N(\nu, \mu),
\end{equation}
\begin{multline}
 \qquad\qquad\qquad
 u \otimes v \mapsto 
 \sum_{(u),(v)} u_{(0)} \otimes v_{(0)} \otimes u_{(1)}v_{(1)} \\ 
 (u \in M(\lambda, \nu), v \in N(\nu, \mu), 
 \lambda, \mu, \nu \in {\V}).
 \qquad
\end{multline}
With this operation, the category $\bold{Com}^f_{\HH}$ of all
finite-dimensional right $\HH$-comodules becomes a 
monoidal (or tensor) category whose unit object $\K \V$
is given by  
\begin{equation}
\label{unitcomdef}		
 \K \V = \bigoplus_{\mu \in \V} \K \mu; 
 \quad
 \mu \mapsto 
 \sum_{\lambda \in \V} \lambda \otimes \emal \enam
 \quad (\mu \in \V).
\end{equation}
If, in addition, $\HH$ is CQT, 
then $\bold{Com}^f_{\HH}$ becomes a
braided monoidal category with braiding
$c_{MN}\!: M \bar{\otimes} N 
\cong N \bar{\otimes} M$ 
given by 
\begin{multline}
\label{cdef}
 \qquad\qquad
 c_{MN} (u \otimes v)
 = \sum_{(u),(v)} v_{(0)} \otimes u_{(0)} {\cal R}^+(u_{(1)},v_{(1)}) 
 \\ (u \in M(\lambda, \nu), v \in N(\nu, \mu), 
 \lambda, \mu, \nu \in {\V}).
 \qquad  
\end{multline}
Next, suppose that $\HH$ has a bijective antipode and 
$M$ is finite-dimensional. 
Then there exists a unique right $\HH$-comodule $M \spcheck$
such that whose underlying vector space is the dual of $M$
and that the coaction satisfies 
\begin{equation}
\label{leftdualdef}
 \sum_{(u)} \langle v, u_{(0)} \rangle S(u_{(1)}) =
 \sum_{(v)} \langle v_{(0)}, u \rangle v_{(1)}
\quad (u \in M, v \in M \spcheck).
\end{equation}
As vector spaces, we have
\begin{equation}
\label{Mcheckij}
 M \spcheck (\lambda, \mu) \cong M(\mu, \lambda)^*
 \quad (\lambda, \mu \in \V).
\end{equation}
The comodule $M \spcheck$ becomes a left dual object of
$M$ via maps
\begin{equation}
\label{bMdef}
 b_M\!:\,
 \K \V \to M \bar{\otimes} M \spcheck ; 
\quad
 \lambda \mapsto 
 \sum_{\mu \in \V}
 \sum_{\p \in \EuScript{M}_{\lambda \mu}} 
 \p \otimes \p \spcheck 
 \quad (\lambda \in \V), 
\end{equation}
\begin{equation}
\label{dMdef}
 d_M\!:\,		
 M \spcheck \bar{\otimes} M \to \K \V\, ;
\quad		 
 \p \spcheck \otimes \q \mapsto
 \delta_{\p \q}\, \mu
 \quad 
 (\p \in \EuScript{M}_{\nu \lambda}, 
 \q \in \EuScript{M}_{\nu \mu}),
\end{equation}
where $\EuScript{M}_{\lambda \mu}$ denotes a 
basis of $M (\lambda, \mu)$ and 
$\{ \p \spcheck\, |\, \p \in 
\EuScript{M}_{\lambda \mu} \}$
its dual basis.
By replacing $S$ in \eqref{leftdualdef}
with $S^{-1}$, we obtain another $\HH$-comodule structure
on $M^*$, which gives the right dual $M^{\land}$ of $M$.
We note that the canonical linear isomorphism 
$I_{M^{\land}}\!: M^{\land} \to M^{\lor}$ satisfies
\begin{equation}
\label{IXu}
 I_{M^{\land}} (X u) =
 (S^{-2})^* (X) I_{M^{\land}} (u) 
\end{equation}
for each $X \in \HH^*$ and $u \in M^{\land}$. 
Here, as usual, we regard $M$ as a left $\HH^*$-module via  
\begin{equation}
\label{Xu}
 X u =
 \sum_{(u)} u_{(0)} \langle X, u_{(1)} \rangle 
 \quad (u \in M, X \in \HH^*).
\end{equation}
Finally, suppose that $\HH$ is a coribbon Hopf $\V$-face algebra.
Then $\bold{Com}^f_{\HH}$ becomes a ribbon category with twist
$\theta_M\!: M \cong M$ given by
\begin{equation}
 \theta_M (u) =
 \cal{V}^{-1} u
 \quad (u \in M).
\end{equation}
\begin{lem}
\label{Trq}					
 Let $\HH$ be a coribbon Hopf $\V$-face algebra 
 with modified ribbon functional $\cal{M}$ such that
 its unit comodule $\K \V$ is absolutely irreducible.
 Then, the quantum trace of the ribbon category
 $\bold{Com}^f_{\HH}$ is given by
 \begin{align}
  \mathrm{Tr}_q (f) & =
  \frac{1}{\mathrm{card} (\V)}  \mathrm{Tr} (\cal{M} f) \\
  & = 
 \label{Trqf}		
 \mathrm{Tr} \left( (\cal{M} f) |_{M(\lambda, -)} \right)
 \end{align}
 for each $M \in \mathrm{ob} ( \bold{Com}^f_{\HH} )$ 
 and $f \in \End_{\HH} (M)$, 
 where $\lambda$ is an arbitrary element of $\V$ and 
 $M ( \lambda, -) = \bigoplus_{\mu} M (\lambda, \mu)$. 
\end{lem}
\begin{pf}
 We define $u_M\!: M \to M^{\lor \lor}$ by
 the composition
 \begin{multline*}
  \qquad
  M \cong M \bar{\otimes} \K \V 
  @>{\id \bar{\otimes} b}>>
  M \bar{\otimes} M^{\lor} \bar{\otimes} M^{\lor \lor} 
  @>{c \bar{\otimes} \id}>> \\
  M^{\lor} \bar{\otimes} M \bar{\otimes} M^{\lor \lor} 
  @>{d \bar{\otimes} \id}>>
  \K \V \bar{\otimes} M^{\lor \lor}  
  \cong M^{\lor \lor}.
  \qquad
 \end{multline*}
 Then, as a consequence of the fact that
 $\bold{Com}_{\HH}^f$ is a rigid braided monoidal category, 
 we have 
 \begin{equation}
 \label{dMcMM}
 d_M \circ c_{M M^{\lor}} =
 d_{M^{\lor}} \circ ( u_M \bar{\otimes} \id_{M^{\lor}} ).
 \end{equation}
 On the other hand, since $u_M (v) = I_M (\cal{U}_1 v)$
 $(v \in M)$, we have
 \begin{equation}
 \label{utheta}
  (u_M \circ \theta_M) (v) = I_M (\cal{M} v)
  \quad (v \in M),
 \end{equation}
 where $I_M$ is as in \eqref{IXu}.
 Hence, by the definition of $\Tr_q$, we obtain 
 \begin{align*}
  \Tr_q (f) (\lambda) & =
  d_{M^{\lor}} \circ ( u_M \theta_M f \bar{\otimes} \id_{M^{\lor}} )
  \circ b_M (\lambda) \\
  = & \sum_{\mu \in \V}
  \sum_{\p \in \EuScript{M}_{\lambda \mu}}
  \langle I_M ( \cal{M} f (\p) ), \p^{\lor} 
  \rangle\, \lambda \\
  & = 	
  \mathrm{Tr} \left( (\cal{M} f) |_{M(\lambda, -)} \right)
  \lambda
 \end{align*}
 as required, 
 where $\EuScript{M}_{\lambda \mu}$ and $\p^{\lor}$ are
 as in \eqref{bMdef}.
\end{pf}		
We say that a right comodule $M$ of 
a $\V$-face algebra $\HH$ is {\it group-like} if 
$\dim (M(\lambda, \mu)) = \delta_{\lambda \mu}$
for each $\lambda, \mu \in \V$.
For each $g \in \mathrm{GLE} (\HH)$, we can define 
a group-like comodule 
$\K \V g = \bigoplus_{\mu \in \V} \K \mu g$				
by the coaction 
$\mu g \mapsto \sum_{\lambda} \lambda g 
\otimes \emal \enam g$.
Conversely, we have the following.
\begin{lem}[\cite{fb}]
\label{glcom}
 For each group-like comodule $M$ and its basis
 $g_{\mu} \in M(\mu, \mu)$ $(\mu \in \V)$,
 we obtain a group-like element
 $g$ by 
 $\emal \enam g = g \binom{\lambda}{\mu}$ and 
 $g_{\mu} \mapsto \sum_{\lambda} g_{\lambda} 
 \otimes g \binom{\lambda}{\mu}$.
 Moreover, every group-like element is obtained in this manner.
\end{lem}
Let $\HH$ be a $\V$-face algebra and 
$\{ L_{\psi}\, |\, \psi \in \LLam \}$ 
its finite-dimensional comodules 
such that $\HH \cong \bigoplus_{\psi \in \LLam} 
\End ( L_{\psi})^*$ as coalgebras.
Let $g$ be a central group-like element of $\HH$.
We say that $g$ is {\it simply reducible} if there exists 
a set $\overline{\LLam}$ and a bijection 
$\varphi\!: \overline{\LLam} \times \Z_{\geq 0}$ 
$\cong$ $\LLam$ such that 
$L_{\varphi (\lambda, n)}$ $\cong$
$\K \V g^n \bar{\otimes} L_{\varphi (\lambda, 0)}$ 
for each $\lambda \in \overline{\LLam}$ and $n \geq 0$.
\begin{lem}[\cite{gal}]
\label{sympredgle}
 Let $\HH$, $g$ etc. be as above. Then, 
 \rom{(1)} the element $g$ is not a zero-divisor of $\HH$.
 \rom{(2)} The quotient $\overline{\HH}= \HH / (g-1)$ 
 is isomorphic to
 $\bigoplus_{\lambda \in \overline{\LLam}} \End (L_{\lambda})^*$
 as coalgebras, where $L_{\lambda}$ is 
 $L_{\varphi (\lambda, 0)}$ viewed as an $\overline{\HH}$-comodule. 
 As $\overline{\HH}$-comodules, we have 
 $L_{\varphi (\lambda, n)} \cong$
 $L_{\lambda}$ for each $\lambda \in \overline{\LLam}$
 and $n$.
\end{lem}			

\section{compact face algebras}


Let $\HH$ be a $\V$-face algebra over the complex
number field $\Bbb C$ and $\times\!: \HH \to \HH;$
$a \mapsto a^{\times}$
an antilinear map.
We say that $(\HH, \times)$ is a {\it costar}
$\V$-{\it face algebra} \cite{cpt} if
\begin{equation}
 (a^{\times})^{\times} = a,
 \quad
 (ab)^{\times} = a^{\times} b^{\times},
\end{equation}
\begin{equation}
 \Delta(a^{\times})
 = \sum_{(a)} a_{(2)}^{\times}
 \otimes a_{(1)}^{\times},
\end{equation}
\begin{equation}
 \enal^{\times} = \emal
\end{equation}
for each $\lambda \in \V$ and $a, b \in \HH$.
Let $[t_{pq}]_{p,q \in I}$ be a finite-size
matrix whose entries are elements of $\HH$.
Then $[t_{pq}]$ is called
a {\it unitary matrix corepresentation}
if $\Delta(t_{pq}) = \sum_{r} t_{pr} \otimes t_{rq}$,
$\cu(t_{pq}) = \delta_{pq}$
and $t_{pq}^{\times} = t_{qp}$.
A costar $\V$-face algebra $\HH$ is called {\it compact}
if $\HH$ is spanned by entries of unitary matrix
corepresentations (cf. \cite{Koornwinder,cpt}).
For each costar $\V$-face algebra $\HH$, its dual 
$\HH^*$ becomes a $*$-algebra via 
\begin{equation}
 \langle X^*, a \rangle = 
 \overline{ \langle X, a^{\times} \rangle}
 \quad (X \in \HH^*, a \in \HH).
\end{equation}
When $\HH$ is a Hopf $\V$-face algebra, 
we also set $a^* = S( a^{\times} )$ and 
$X^{\times} = S^* (X)^*$ for each $a \in \HH$
and $X \in \HH^*$.
These operations satisfy
\begin{equation}
\label{a**}
 (a^*)^* = a,
\quad
 (X^{\times})^{\times} = X
 \quad (a \in \HH, X \in \HH^*).
\end{equation}
For each compact Hopf $\V$-face algebra $\HH$,
its {\it Woronowicz functional} \cite{Koornwinder,cpt}
$\cal{Q} \in \HH^*$ is a unique group-like
functional such that
\begin{equation}
\label{QXQ-}
 \qquad
 \cal{Q} X \cal{Q}^{-1} =
 (S^2)^* (X)
 \quad (X \in \HH^*),
\end{equation}
\begin{equation}
\label{TrQ}
 \sum_p \cal{Q} (t_{pp}) =
 \sum_p \cal{Q}^{-1} (t_{pp}),
\end{equation}
and that $[\cal{Q} (t_{pq})]$ is a positive matrix, 
where $[t_{pq}]$ denotes an arbitrary
unitary matrix corepresentation.
The functional $\cal{Q}$ also satisfies
\begin{equation}
\label{Q*}
 \cal{Q}^* =  \cal{Q},
\quad
 \cal{Q}^{\times} =  \cal{Q}^{-1}.
\end{equation}
Let $(\HH,\times)$ be a costar $\V$-face algebra
and $M$ a finite-dimensional $\HH$-comodule
equipped with a Hilbert space structure $(\,|\,)$.
We say that $(M,(\,|\,))$ is {\it unitary} if
\begin{equation}
 \qquad
 \sum_{(u)} (u_{(0)}\, |\, v)\, u_{(1)} =
 \sum_{(v)} (u\, |\, v_{(0)})\, v_{(1)}^{\times}
 \quad (u, v\in M).
\end{equation}
We note that $(M,(\,|\,))$ is unitary 
if and only if \eqref{Xu} gives a $*$-representation 
of $\HH^*$.
A costar $\V$-face algebra $(\HH,\times)$ is compact
if and only if every finite-dimensional $\HH$-comodule
is unitary for some $(\,|\,)$.

\begin{prop} Let $\HH$ be a compact $\V$-face algebra.		
 \rom{(1)} The unit comodule $\C \V$ becomes a unitary comodule
 via $( \lambda\, |\, \mu ) = \delta_{\lambda \mu}$. 
 \rom{(2)} For each unitary comodules		
 $M$ and $N$, $M \bar{\otimes} N$ becomes a unitary comodule
 via $(u \otimes v\, |\, u^{\prime} \otimes v^{\prime} )$ $=$ 
 $(u\, |\, u^{\prime}) (v\, |\, v^{\prime} )$
 $(u \in M(\lambda, \nu), v \in N(\nu, \mu),   
  u^{\prime} \in M(\lambda^{\prime}, \nu^{\prime}), 
  v^{\prime} \in N(\nu^{\prime}, \mu^{\prime}) )$.
 \rom{(3)} If $\HH$ has an antipode, then the left dual 
 $M \spcheck$ becomes a unitary comodule via 
 \begin{equation}
  (u\, |\, v) = 
  ( \cal{Q} \Upsilon^{-1} (v)\, |\, \Upsilon^{-1} (u) ),
 \end{equation}					
 where $\Upsilon\!: M \to M \spcheck$ 
 denotes the antilinear isomorphism
 defined by 
 \begin{equation}
  \langle \Upsilon (u),\, v \rangle 
  = (v\, |\, u) 
 \quad (u, v \in M).
 \end{equation}					
\end{prop}
\begin{pf}
 The proof of
 Part (1) and Part (2) is straightforward.
 Using the second equality of \eqref{a**}, 
 we obtain $(S^2)^* (X^{\times}) =$
 $S^* (X^*)$. Using this together with
 \eqref{QXQ-} and 
 \begin{equation}
 \label{XUu}
 X \Upsilon (u) = \Upsilon (X^{\times} u)
 \quad  (u \in M, X \in \HH^*),
 \end{equation}
 we obtain
 \begin{align*}
  (u\, |\, X v) 
  & = 
  ( S^* (X^*) \cal{Q} \Upsilon^{-1} (v)\, |\, \Upsilon^{-1} (u) ) \\
  & = 
  ( \cal{Q} \Upsilon^{-1} (v)\, |\, (X^*)^{\times} \Upsilon^{-1} (u) ) \\
  & = (X^* u\, |\, v)
 \end{align*}					
 for each $u, v \in M \spcheck$, as required.
\end{pf}

By the proposition  above, the category $\bold{Com}_{\HH}^{fu}$
of all finite-dimensional unitary comodules of 
a compact (Hopf) $\V$-face algebra $\HH$ becomes
a (rigid) monoidal category.
Let $M$ and $N$ be unitary $\HH$-comodules 
and $f\!: M \to N$
an $\HH$-comodule map. We define a linear map
$\bar{f}\!: N \to M$ by
$(\overline{f} (n)\, |\, m) = (n\, |\, f(m))$.
We have
\begin{equation}
 \overline{\overline{f}} = f,
\quad
 \overline{f \circ g} =
 \overline{g} \circ \overline{f},
\quad
 \overline{f \bar{\otimes} g} =
 \overline{f} \bar{\otimes} \overline{g},
\end{equation}
\begin{equation}
 \overline{f + g} =
 \overline{f} + \overline{g}, 
\quad
 \overline{c f} = \overline{c} \overline{f}
\quad (c \in \C).
\end{equation}
\begin{prop}
For each compact Hopf $\V$-face algebra $\HH$ and 
its comodule $M$, we have
\begin{equation}
\label{dMbar}
 \overline{d_M} = \left( 
 (I_{M^{\land}} \circ \cal{Q}) \bar{\otimes} \id_M \right) 
 \circ b_{M^{\land}},
\end{equation}
\begin{equation}
\label{bMbar}
 \overline{b_M} = d_{M^{\land}} \circ
 ( \id_M \bar{\otimes} ( (I_{M^{\land}} \circ \cal{Q})^{-1} ).
\end{equation}
\end{prop}		
\begin{pf}	
 Using \eqref{Q*} and \eqref{XUu}, we obtain
 \begin{equation}
 \label{q|Qp}
  ( \q^{\lor}\, |\, \cal{Q} \p^{\lor}) 
   = 
  \overline{ \langle \p^{\lor},\, 
  \Upsilon^{-1} (\q^{\lor}) \rangle},
 \quad
 ( \Upsilon (\q)\, |\, \r^{\lor}) 
   = 
  \overline{ \langle \r^{\lor},\, 
  \cal{Q} \q \rangle},
 \end{equation}					
 where $\{ \p \}$ and $\{ \p^{\lor} \}$ is as in 
 \eqref{bMdef}. 
 Using the first equality of  \eqref{q|Qp}, we obtain
 \begin{align*}
  \left(	
  \q^{\lor} \otimes \r\, |\, 
  ( (I_{M^{\land}} \circ \cal{Q}) \bar{\otimes} \id_M ) 
  \circ b_{M^{\land}} (\lambda) \right) 
  & = 
  \sum_{\mu \in \V} \sum_{\p \in \EuScript{M}_{\mu \lambda}}
  \left(					
  \q^{\lor} \otimes \r\, |\, 
  \cal{Q} \p^{\lor} \otimes \p 
  \right)
  \\ = 
  \sum_{\mu \in \V}
  \sum_{\p \in \EuScript{M}_{\mu \lambda}}
  \left( \r\, |\, 
  \langle \p^{\lor}, \Upsilon^{-1} (\q^{\lor}) 
  \rangle \p \right) 
  & = 
  \, \delta_{\lambda \en (\q)}\, \delta_{\q \r}, 
 \end{align*}					
 where the first equality follows from 
 \eqref{IXu}.
 This proves \eqref{dMbar}.
 Similarly, \eqref{bMbar} follows from the second equality
 of \eqref{q|Qp}.
\end{pf}
Let $\HH$ be a costar CQT $\V$-face algebra.
We say that $\HH$ is of {\it unitary type}
if
\begin{equation}
 \label{R(atbt)}
 \Rp (a^{\times}, b^{\times}) =
 \overline{\Rm (a, b)} 
\quad (a, b \in \HH).
\end{equation}
In this case, the Drinfeld functionals of $\HH$ satisfy
\begin{equation}
\label{U*}
 \cal{U}_1^* =  \cal{U}_2, 
\quad
 \cal{U}_2^* =  \cal{U}_1.
\end{equation}
\begin{rem}
 When $q > 0$, the function algebras of the usual 
 quantum groups (such as $\mathrm{Fun}(SL_q (N))$,
 $\mathrm{Fun}(Sp_q (2N))$) are both 
 CQT and compact. However, they are not of unitary type 
 but rather of ``Hermitian type.''
\end{rem}
\begin{prop}[\cite{cpcq}]
\label{Q=M}
 For each compact CQT Hopf $\V$-face algebra $\HH$ of 
 unitary type, its Woronowicz functional is a modified 
 ribbon functional of $\HH$. 
 Moreover, the corresponding ribbon functional
 $\cal{V}_{\cal{Q}}$ satisfies
 \begin{equation}
 \label{V*}	
  \cal{V}_{\cal{Q}}^{\;*} = \cal{V}_{\cal{Q}}^{-1}.
 \end{equation}		
 %
\end{prop}
We call $\cal{V}_{\cal{Q}}$ 
the {\it canonical ribbon functional} of $(\HH, \times)$. 
We note that the expression 
$\cal{U}_1 = \cal{V}_{\cal{Q}} \cal{Q}$
gives the ``polar decomposition'' of $\cal{U}_1$.

Let $\HH$ be a compact CQT Hopf $\V$-face algebra of unitary type,
equipped with the canonical ribbon functional.
By \eqref{R(atbt)} and \eqref{V*}, 
we have
\begin{equation}
 \overline{c_{M N}} = (c_{M N})^{-1},
\quad
 \overline{\theta_M} = \theta_M^{-1}.
\end{equation}
Moreover, by \eqref{dMcMM}, \eqref{utheta}, 
\eqref{bMbar}, we have
\begin{equation}
 \overline{b_M} = d_M \circ c_{M M^{\lor}} 
 \circ
 ( \theta_M \bar{\otimes} \id_{M^{\lor}} ). 
\end{equation}
Furthermore, by \eqref{Trqf} and the defining properties
of $\cal{Q}$, we have 
\begin{equation}
 \Tr_q (f \overline{f}) =
 \frac{1}{\card (\V)} 
 \Tr ( \overline{f} \cal{Q} f ) > 0
\end{equation}
for every $f \ne 0$.
Thus, the category $\bold{Com}_{\HH}^{fu}$ 
is a unitary ribbon category (cf. \cite{Turaev}).

\begin{exam}
 Let $L=1$ and $\SS$ as in Example \ref{L=1}. Then
 $\SS$ is a compact CQT Hopf $\V$-face algebra of unitary type
 with costar structure
 $e^i_j (m)^{\times} = e^j_i (m)$.
 Its Woronowicz functional agrees with
 the counit $\cu$.
\end{exam}


\section{flat face models}

					
Let $\G$ be a finite oriented graph with set of vertices $\V$.
We say that a quadruple 
$\left( \r \frac[0pt]{\p}{\q} \s \right)$ 
or a diagram
\begin{equation}
 \begin{CD}
  \lambda @>{\p}>> \mu \\
  @V{\r}VV @VV{\s}V \\
  \nu @>{\q}>> \xi \\
 \end{CD}
 \label{face}
\end{equation}
is a {\it face} if 
$\p,\q, \r, \s \in {\G}^1$ and
\begin{equation} 
 \st (\p) = \lambda = \st (\r),
 \quad \en (\p) = \mu = \st (\s),  
 \quad \en (\r) = \nu = \st (\q), 
 \quad \en (\q) = \xi = \en (\s).
 \label{facecond}
\end{equation}
When $\G$ has no multiple edge, we also write
$\face{\r}{\p}{\q}{\s} =$
$\left( {\lambda \, \mu}
\atop {\nu \, \xi} \right)$.
We say that $(\G, w)$ is a ($\V$-){\it face model} 
over a field $\K$ if 
$w$ is a map which assigns a number
$w \!\! \left[ \r \frac[0pt]{\p}{\q} \s \right]$
$\in \K$ to each face 
$\left( \r \frac[0pt]{\p}{\q} \s \right)$
of $\G$.
We call $w$ {\it Boltzmann weight} of $(\G, w)$.
For convenience, we set 
$w \!\! \left[ \r \frac[0pt]{\p}{\q} \s \right] = 0$
unless
$\left( \r \frac[0pt]{\p}{\q} \s \right)$
is a face.   
For a face model $(\G,w)$, we identify $w$ 
with the linear operator on
$\K \G^2$ $=$ $\bigoplus_{\p \in \G^2} \K \p$ 
given by
\begin{equation}
 w(\p \cdot \q) =
 \sum_{\r \cdot \s \in {\G}^2}
 w \!\! \left[ \r \frac[0pt]{\p}{\s} \q \right]
 \r \cdot \s 
 \quad (\p \cdot \q \in \G^2).
\end{equation}
For $m \geq 2$ and $1 \leq i < m$, we define an 
operator $w_i = w_{i/m}$ on $\K \G^m$ by
$w_{i/m} ( \p \cdot \q \cdot \r )$
$=$ $\p \otimes w (\q) \otimes \r$ 
$(\p \in \G^{i-1}, \q \in \G^{2}, \r \in \G^{m-i-1})$, 
where we identify $(\p_1, \ldots, \p_m) \in \G^m$
with 
$\p_1 \otimes \ldots \otimes \p_m \in 
(\K \G^1)^{\otimes m}$.

A face model is called {\it invertible}
if $w$ is invertible as an operator on  
$\K \G^2$.
An invertible face model is called {\it star-triangular}
(or {\it Yang-Baxter}) if $w$ satisfies
the braid relation $w_1 w_2 w_1$ $=$
$w_2 w_1 w_2$ in $\End (\K \G^3)$.

%
			
For a star-triangular face model $(\G, w)$, 
the operators $w_{i/m}$ $(1 \leq i < m)$ define
an action of 
the $m$-string braid group $\frak{B}_m$ on $\K \G^m_{\lambda \mu}$
for each $m \geq 2$ and $\lambda, \mu \in \V$. 					
The following proposition gives a face version of the 
FRT construction.
\begin{prop}[\cite{RTF,LarsonTowber,gd,fb}]
\label{FRT}
Let $(\G,w)$ be a $\V$-face model and $\HH (\G)$ as in $\S$ 1.
Let ${\frak I}$ be an ideal of $\HH (\G)$ generated by the following 
elements\rom:
\begin{equation}
 \sum_{\r \cdot \s \in {\G}^2}
 w \!\! \left[ \r \frac[0pt]{\p}{\s} \q \right]
 e{\aaa \cdot \bbb  \choose \r \cdot \s} -
 \sum_{\ccc \cdot \ddd \in {\G}^2}
 w \!\! \left[ \aaa \frac[0pt]{\ccc}{\bbb} \ddd \right] 
 e{\ccc \cdot \ddd  \choose \p \cdot \q}
 \quad (\p \cdot \q,\, \aaa \cdot \bbb \in \G^2).
\end{equation}
Then ${\frak I}$ is a coideal of $\HH (\G)$ and the quotient
$\Aw := \HH (\G) / {\frak I}$ 
becomes a $\V$-face algebra.   
If $(\G,w)$ is star-triangular, then there exist unique 
bilinear pairings $\Rpm$ on $\Aw$
such that $(\Aw, \Rpm)$ is a CQT $\V$-face algebra
and that 
\begin{equation}
\label{UnivR=w}
  \Rp \left( e{\p \choose \q}, e{\r \choose \s} \right)
 =  w \!\! \left[ \r \frac[0pt]{\q}{\p} \s \right]
\end{equation}
for each $\p,\q, \r, \s \in {\G}^1$. 
\end{prop}

We denote the image of $e \binom{\p}{\q}$ by the projection
$\HH (\G) \to \Aw$ again by $e \binom{\p}{\q}$.
Then $\frak{A}^m (w)$ $:=$
$\sum_{\p, \q \in \G^m} \K e \binom{\p}{\q}$ 
becomes a subcoalgebra of $\Aw$ for each $m \geq 0$.			
As the usual FRT construction (cf. \cite[Proposition 2.1]{qcg}),  
we have the following.
\begin{prop}
\label{Aw*}
 For each star-triangular face model $(\G, w)$,
 we have $\frak{A}^m (w)^* \cong \Hom_{\frak{B}_m} (\K \G^m)$
 $(m \geq 2)$ as $\K$-algebras.
\end{prop}		

We say that $\left( \r \frac[0pt]{\p}{\q} \s \right)$
or \eqref{face} is a {\it boundary condition} of size $m \times n$
if $\p,\q \in {\G}^n, \r, \s \in {\G}^m$ 
and the relation \eqref{facecond}
is satisfied for some $\lambda, \mu, \nu, \xi$.
For a face model $(\G,w)$, 
we define its {\it partition function}
to be an extension
$w: \{ \text{ boundary conditions of size} 
\; m \times n ; \; m,n \geq 1 \} \to \K$ 
of the map $w$ which is determined by the following two recursion relations:

\begin{gather}
 w \!\! \left[ \r \frac[0pt]{\p \cdot {\p}^{\prime}}{\q \cdot {\q}^{\prime}} \s \right]
 =
 \sum_{\aaa \in {\G}^m}
 w \!\! \left[ \r \frac[0pt]{\p}{\q} \aaa \right] 
 w \!\! \left[ \aaa \frac[0pt]{{\p}^{\prime}}{{\q}^{\prime}} \s \right], \\
 w \!\! \left[ \r \cdot {\r}^{\prime} \; \frac[0pt]{\p}{\q} 
 \; \s \cdot {\s}^{\prime} \right]
 =
 \sum_{\aaa \in {\G}^n}
 w \!\! \left[ \r \frac[0pt]{\p}{\aaa} \s \right] 
 w \!\! \left[ {\r}^{\prime} \frac[0pt]{\aaa}{\q} {\s}^{\prime} \right] \\
 \qquad\qquad\qquad\qquad\qquad\qquad\quad
 \left( \p,\q \in {\G}^n, {\p}^{\prime},{\q}^{\prime} 
 \in {\G}^{n^{\prime}},
 \r,\s \in {\G}^m, {\r}^{\prime},{\s}^{\prime} 
 \in {\G}^{m^{\prime}} \right). \nonumber
\end{gather}
Also, we set 
$w \!\! \left[ \r \frac[0pt]{\p}{\q} \s \right] = {\delta}_{\p\q}$
(respectively
$w \!\! \left[ \r \frac[0pt]{\p}{\q} \s \right] = {\delta}_{\r\s}$)
if $\r,\s \in {\G}^0$ (respectively $\p,\q \in {\G}^0$). 
With this notation, the relation \eqref{UnivR=w} holds 
for every star-triangular face model $(\G,w)$ and 
$\p,\q \in {\G}^n, \r, \s \in {\G}^m (m,n \geq 0)$. 
Next, we recall the notion of flat face model \cite{gal}, 
which is a variant of A. Ocneanu's notion of 
flat biunitary connection (cf. \cite{Ocneanu}).
Let $(\G, w)$ be an invertible face model 
with a fixed vertex $* \in \V = {\G}^0$. 
We assume that $\LLam_{\G}^m \ne \emptyset$
for each $m \geq 0$ and that 
$\V = \bigcup_{m \geq 0} \V (m)$, where
$\LLam_{\G}^m$ $=$ $\LLam_{\G *}^m$ and
$\V (m)$ $=$ $\V (m)_{\G *}$ are defined by  
\begin{align}
\label{LambdaGdef}
 & \LLam_{\G} = 
 \bigl\{ (\lambda, m) \in \V \times \Z_{\geq 0} \bigm| 
 \G^m_{* \lambda} \ne \emptyset \bigr\}, \\
 & \LLam_{\G}^m = 
 \LLam_{\G} \cap ( \V \times \{ m \} ), \\
 & \V (m) = 
 \bigl\{ \lambda \in \V \bigm|
 (\lambda, m) \in \LLam_{\G} \bigr\}.
\end{align}
For each $m \geq 0$, we define the algebra 
$\Str{m}$ by
\begin{equation}
 \Str{m} = 
 \prod_{\lambda \in \V (m)}
 \End(\K \G^m_{* \lambda})
\end{equation}
and call it {\it string algebra} of $(\G, w, *)$.
For each $m, n \geq 0$, we define the algebra map
$\iota_{mn}\!: \Str{m} \to \Str{m+n}$ by 
$\iota_{mn} (x) (\p \cdot \q) = x\p \otimes \q$
$(\p \in \G^m_{* \lambda}, \q \in \G^n_{\lambda \mu})$. 
For each $1 \leq i < m$, we define the element
$\wsta_i = \wsta_{i/m}$ of $\Str{m}$ to be 
the restriction of $w_{i/m}$ on $\K \G^m_{*-}$,
where $\G^m_{*-}$ $=$ 
$\coprod_{\lambda \in \V} \G^m_{* \lambda}$.
We say that $(\G, w, *)$ is a {\it flat face model} if 
\begin{equation}
\label{flatdef}
 \iota_{mn} (x) \wsta_{nm} \iota_{nm} (y) \wsta_{nm}^{-1} =
 \wsta_{nm} \iota_{nm} (y) \wsta_{nm}^{-1}  \iota_{mn} (x)
\end{equation}
for each $x \in \Str{m}$ and $y \in \Str{n}$
$(m, n \geq 0)$, 
where $\wsta_{mn} \in \Str{m+n}$ is defined by
\begin{equation}
\label{wstamndef}
 \wsta_{mn} =  
 (\wsta_n \wsta_{n+1} \cdots \wsta_{m+n-1})
 (\wsta_{n-1} \wsta_{n} \cdots \wsta_{m+n-2})  
 \cdots
 (\wsta_{1} \wsta_{2} \cdots \wsta_{m}).  
 \end{equation}
For each flat $\V$-face model $(\G, w, *)$,
$n \geq 0$ and $\lambda, \mu \in \V$, 
there exists a unique left action $\Gamma$ of 
$\Str{n}$ on $\K \G^n_{\lambda \mu}$ such that
\begin{equation}
\label{Gammadef}
 \p \otimes (\Gamma (x) \q) =
 \wsta_{nm} \iota_{nm} (x) \wsta_{nm}^{-1} 
 (\p \cdot \q) 
\end{equation}
for each $m \geq 0$,
$\p \in \G^m_{* \lambda}, \q \in \G^n_{\lambda \mu}$
and $x \in \Str{n}$.
Using this action, we define {\it costring algebra}
\begin{equation}
 \mathrm{Cost} (w, *) =
 \bigoplus_{m \geq 0} \mathrm{Cost}^m (w, *)
\end{equation}
to be the quotient $\V$-face algebra of 
$\bigoplus_{m \geq 0}
\End_{\K} (\K \G^m)^* \cong \HH (\G)$
given by
\begin{equation}
 \mathrm{Cost}^m (w, *) = 
 \End_{\Str{m}} (\K \G^m)^* .
\end{equation}

For each $\lambda, \mu \in \V$ and 
$(\nu, m) \in \LLam_{\G}$,
we define the non-negative integer $N^{\mu}_{\lambda \nu} (m)$
by the irreducible decomposition of $\K \G^m_{\lambda \mu}$:
\begin{equation}
\label{frdef}
 [\K \G^m_{\lambda \mu}] =
 \sum_{\nu \in \V (m)} N^{\mu}_{\lambda \nu} (m)
 [\K \G^m_{* \nu}], 
\end{equation}
where, for each $\Str{m}$-module $V$, 
$[V]$ denotes the element of the 
Grothendieck group $K_0 (\Str{m})$ 
corresponding to $V$ (see e.g. \cite[\S 5.1]{ChariPressley}).
We call $N^{\mu}_{\lambda \nu} (m)$  
{\it fusion rules} of $(\G, w, *)$, 		
\begin{thm}[\cite{gal}]
\label{Coststrthm} 
 Let $(\G, w, *)$ be a flat $\V$-face model with fusion rule 
 $N_{\lambda \nu}^{\mu} (m)$.
 \rom{(1)} For each $(\lambda, m) \in \LLam_{\G}$, up to isomorphism
 there exists a unique right 
 $\mathrm{Cost}^m (w, *)$-comodule $L_{(\lambda, m)}$
 such that 
 \begin{equation} 
 \label{dimLlm*m}
  \dim L_{(\lambda, m)}(*, \mu) = \delta_{\lambda \mu}
 \end{equation} 
 for each $\mu \in \V$.
 As coalgebras, we have
 \begin{equation} 
  \mathrm{Cost}^m (w, *) \cong 
  \bigoplus_{\lambda \in \V (m)}
  \mathrm{End} (L_{(\lambda, m)})^*.
 \end{equation} 
 \rom{(2)} In the corepresentation ring 
 $K_0 (\bold{Com}^f_{\mathrm{Cost} (w, *)})$, we have
 \begin{equation} 
  [L_{(*,0)}] = 1,
 \end{equation}
 \begin{equation} 
  [L_{(\lambda,m)}] [L_{(\mu,n)}] 
  = \sum_{\nu \in \V (m+n)} N^{\nu}_{\lambda \mu} (n) [L_{(\nu, m+n)}].
 \end{equation} 
 Moreover, for each $\mathrm{Cost}^m (w, *)$-comodule $M$, 
 we have
 \begin{equation} 
 \label{decompform}
  [M] = \sum_{\lambda \in \V (m)}
  \dim \left(  M (*, \lambda) \right) [L_{(\lambda, m)}].
 \end{equation} 	
 \rom{(3)} We have 
 \begin{equation} 
 \label{dimLnmlm}
  \dim \left( L_{(\nu, m)} (\lambda, \mu) \right) = 
  N^{\mu}_{\lambda \nu} (m). 
 \end{equation} 
\end{thm}

\begin{lem}
 For each flat star-triangular face model, we have
 \begin{equation}
  \Gamma (\wsta_{i/n}) = w_{i/n}
  \quad (n \geq 2,\, 1 \leq i < n).  
 \end{equation}
\end{lem}
\begin{pf}
 By the braid relation, we have 
 \begin{equation}
 \label{wnmwiwnm}
  \wsta_{nm} \wsta_{i/n+m} \wsta_{nm}^{-1} =
  \wsta_{i+m / n+m} 
\end{equation}
for each $1 \leq i < n$.
Using this together with 
$\iota_{nm} (\wsta_{i/n})$ $=$
$\wsta_{i/n+m}$, 
we obtain
\begin{equation}
\label{pGwq}
 \p \otimes (\Gamma (\wsta_{i/n}) \q) =
 \p \otimes w_{i/n} \q
\end{equation}
for each $m \geq 0$, $\p \in \G^m_{* \lambda}$ and 
$\q \in \G^n_{\lambda \mu}$ as required.
\end{pf}
\begin{prop}
\label{flatcond}
 Let $(\G, w)$ be a star-triangular $\V$-face model 
 with a fixed vertex $* \in \V$.
 Then $(\G, w, *)$ is flat if $\K \G^m_{* \lambda}$ 
 is an absolutely irreducible $\frak{B}_m$-module
 for each $(\lambda, m) \in \LLam_{\G}$. 
 In this case, we have $\mathrm{Cost} (w, *)$
 $=$ $\Aw$ as quotients of $\HH (\G)$.
\end{prop}
\begin{pf}
 Using \eqref{wnmwiwnm}, we see that \eqref{flatdef}
 holds for every $x \in \Str{m}$ and $y = \wsta_i$
 $(1 \leq i < n)$. Hence $(\G, w, *)$ is flat if 
 $\Str{m} = \langle \wsta_1, \ldots, 
 \wsta_{m-1} \rangle$
 for each $m > 1$.
 The second assertion follows from Proposition \ref{Aw*}
 and the lemma above.
\end{pf}


%
%

\section{$SU(N)_L$-SOS models}

%
%
%
In order to construct the algebras $\SS$,
we first recall ${SU(N)}_L$-SOS models 
(without spectral parameter) \cite{JMO},
which are equivalent to H. Wenzl's representations of 
Iwahori-Hecke algebras (cf. \cite{Wenzl}) and also,
the monodromy representations of the braid group
arising from conformal field theory 
(cf. A. Tsuchiya and Y. Kanie \cite{TK}). 
Let $N \geq 2$ and $L \geq 2$ be integers. 
For each $1 \leq i \leq N$, we define the vector
$\hat{i} \in \Bbb{R}^N$ by
$\hat{1} = (1 - 1/N, -1/N, \ldots, -1/N)$,
$\ldots$,
$\hat{N} = (- 1/N, \ldots,  -1/N, 1 -1/N)$.
Let $\V = \V_{NL}$ be the subset of $\Bbb{R}^N$ given by

\begin{equation}
\label{Vdef}
 {\V}_{NL} = 
 \bigl\{ \lambda_1 \hat{1} + \cdots \lambda_N \hat{N} 
 \bigm| \lambda_1, \ldots, \lambda_N \in \Z, 
 L \geq {\lambda}_1 \geq \dots \geq {\lambda}_{N} = 0
 \bigr\}.
\end{equation} 
For $\lambda \in \V$, we define integers
$\lambda_1, \ldots, \lambda_N$ and $| \lambda |$ by
$\lambda = \sum_i \lambda_i \hat{i}$, 
$\lambda_N = 0$ and $| \lambda | = \sum_i \lambda_i$.  
For $m \geq 0$, we define a subset ${\G}^m$ of ${\V}^{m+1}$ by 
\begin{equation}
 {\G}^m = 
 {\V}^{m+1} \cap 
 \bigl\{ \p = (\lambda\, |\, i_1, \ldots, i_m)
 \bigm| \lambda\in \V,\, 
 1 \leq i_1, \ldots, i_m \leq N \bigr\},   
\end{equation}
where for $\lambda \in \Bbb{R}^N$ and 
$1 \leq i_1, \ldots, i_m \leq N$,
we set 
\begin{equation}  
 (\lambda\, |\, i_1, \ldots, i_m) =
 (\lambda, \lambda + \hat{i}_1, \ldots, 
 \lambda + \hat{i}_1 + \cdots + \hat{i}_m).
\end{equation}
Then $(\V,{\G}^1)$ defines an oriented graph $\G = {\G}_{N,L}$ 
and ${\G}^m$ is identified with the set of paths of $\G$ of length $m$.    
For $\p = (\lambda\, |\, i,j)$, 
we set 
${\p}^{\dag} = (\lambda\, |\, j,i)$ and
\begin{equation}
 d(\p)
 = {\lambda}_i - {\lambda}_j + j - i. 
\end{equation} 
We define subsets
${\G}^2[\to]$, ${\G}^2[\;\downarrow\;]$ and
${\G}^2[\searrow]$ of ${\G}^2$ by				
\begin{align}
 {\G}^2[\to] & =
 \bigl\{ \p \in {\G}^2 
 \bigm| \p^{\dag} = \p \bigr\},  
 \\
 {\G}^2[\;\downarrow\;] & =
 \bigl\{ \p \in {\G}^2 
 \bigm| 
 \p^{\dag} \not\in \G^2 \bigr\},
 \\
 {\G}^2[\searrow] & =  
 \bigl\{ \p \in {\G}^2 
 \bigm| \p \not= \p^{\dag} \in \G^2 
 \bigr\}.  
\end{align} 
Let $t \in\C$ be a primitive $2(N+L)$-th root of $1$.
Let $\epsilon$ be either $1$ or $-1$
and $\zeta$ a nonzero parameter.
We define a face model 
$(\G,w_{N,t,\epsilon}) = ({\G}_{N,L},w_{N,t, \epsilon, \zeta})$ by setting
\begin{equation}
 w_{N,t,\epsilon} \!
 \begin{bmatrix}
   \lambda           & \lambda + \hat{i} \\
   \lambda + \hat{i} & \lambda + \hat{i} + \hat{j}
 \end{bmatrix}  
 = - \zeta^{-1} t^{-d(\p)} \frac{1}{[d(\p)]},
\end{equation} 
\begin{equation}
 w_{N,t,\epsilon} \!
 \begin{bmatrix}
   \lambda           & \lambda + \hat{i} \\
   \lambda + \hat{j} & \lambda + \hat{i} + \hat{j}
 \end{bmatrix}  
 = \zeta^{-1} \epsilon \, \frac{[d(\p)-1]}{[d(\p)]}, 
\end{equation}
\begin{equation}
 w_{N,t,\epsilon} \!
 \begin{bmatrix}
   \lambda           & \lambda + \hat{k} \\
   \lambda + \hat{k} & \lambda + 2 \hat{k}
 \end{bmatrix}  
 = \zeta^{-1} t
\end{equation}
for each 
$\p = (\lambda\, |\, i,j) \in 
\G^2[\searrow] \amalg \G^2 [\,\downarrow\,]$
and 
$(\lambda\, |\, k,k) \in \G^2[\to]$, 
where $[n] = (t^n-t^{-n})/(t-t^{-1})$
for each $n \in \Z$.
We call $(\G,w_{N,t,\epsilon})$ {\it $SU(N)_L$-SOS model} 
(without spectral parameter) \cite{JMO}.
It is known that $(\G,w_{N,t,\epsilon})$ is star-triangular.
Moreover, H. Wenzl \cite{Wenzl} showed that 
$\C \G^m_{0 \lambda}$ is an irreducible 
$\frak{B}_m$-module for each $m \geq 0$ and 
$\lambda \in \Lambda^m_{\G 0}$. 
Therefore $(\G,w_{N,t,\epsilon},0)$ is flat 
by Proposition \ref{flatcond}.
In \cite{GoodmanWenzl}, F. Goodman and H. Wenzl showed that 
the fusion rule of $(\G,w_{N,t,\epsilon},0)$ agrees with that of
$SU(N)_L$-WZW model. We give another proof of 
their result in the next section.
\begin{rem}
\label{rembase}
(1)
Strictly speaking, Wenzl deals with $w_{N,t,\epsilon}$
only when $\epsilon = - 1$.
However, it is clear that his arguments are applicable to
the case $\epsilon = 1$.				
The results for $\frak{A}(w_{N,t,1})$ also follows from
those of $\frak{A}(w_{N,t,-1})$, since the former is a 
2-cocycle deformation of the latter (cf. \cite{DT}).
Hence, $\epsilon$ may be viewed as a gauge parameter. 
\\			
(2)
In order to avoid using square roots of complex numbers, 
we use a different normalization of $w_{N,t, -1}$
from Wenzl \cite{Wenzl}.
For each $\p \in \G^m$ $(m \geq 1)$, 
we define $\kappa (\p) \in \C$ by
\begin{gather}
\label{kappadef}
 \kappa (\p \cdot \q) = \kappa (\p) \kappa (\q) 
 \quad (\p \in \G^m, \q \in \G^n, m, n > 0), \\
 \kappa (\lambda\, |\, i) = 			
 \prod_{k = i+1}^{N} A_{d(0\, |\, i,k)+1}
 A_{d(0\, |\, i,k)+2} \cdots A_{d(\lambda\, |\, i,k)}
 \quad ((\lambda\, |\, i) \in \G^1),
 \nonumber 
\end{gather}
where $A_d = a_d / \sqrt{a_d a_{-d}}$ and 
$a_d = [d+1] / [2] [d]$.
Note that $\kappa (\p)$ satisfies
\begin{equation}
 \frac{\kappa (\p)}{\kappa (\p^{\dag})} =
 A_{d (\p)}
\end{equation}
for each $\p \in \G^2 [\searrow]$.
By replacing the basis $\{ \p \}$ of $\C \G^m$
with $\{ \kappa (\p) \p \}$,
we obtain Wenzl's 
original expression
of the Hecke algebra representation.
It is also useful to use $\{ \kappa (\p)^2 \p \}$ instead of
$\{ \p \}$ (see \S 12). The corresponding Boltzmann weight 
$w^{\Sigma}_{N,t,\epsilon}$ satisfies
\begin{equation}
 w^{\Sigma}_{N,t,\epsilon} \!
 \left[ \r \frac[0pt]{\p}{\s} \q \right] :=
 \left( 
 \frac{\kappa (\p \cdot \q)}{\kappa (\r \cdot \s)} 
 \right)^2 
 w_{N,t,\epsilon} \!
 \left[ \r \frac[0pt]{\p}{\s} \q \right] =
 w_{N,t,\epsilon} \!
 \left[ \p \frac[0pt]{\r}{\q} \s \right].
\end{equation} 
We call $\{ \p \}$ and $\{ \kappa (\p)^2 \p \}$ 
{\it rational basis of type} $\Omega$ 
and {\it type} $\Sigma$ respectively.   
\end{rem}

\section{The algebra $\SS$}
Applying Proposition \ref{FRT} to $(\G,w_{N,t,\epsilon, \zeta})$, 
we obtain a CQT $\V$-face algebra $\ASe$ $=$
$\mathrm{Cost} (w_{N,t,\epsilon}, 0)$.
In order to define the ``(quantum) determinant'' of $\AS$,
we introduce an algebra $\Omega = \Omega_{N,L, \epsilon}$,
which is a face-analogue of the exterior algebra.
It is defined by generators $\omega (\p)$ 
$(\p \in \G^m; m \geq 0)$
with defining relations: 
\begin{equation}
 \sumk \omega (k) \, = \, 1,
\end{equation}
\begin{equation}
 \omega (\p) \omega (\q) \, = \,
 \delta_{\en (\p) \st (\q)} \, \omega (\p \cdot \q),
\end{equation}
\begin{equation}
 \omega (\p) = - \epsilon\, \omega (\p^{\dag})
 \quad (\p \in \G^2 [\searrow]),
\end{equation}
\begin{equation}
 \omega (\p)  = 0
 \quad (\p \in \G^2[\to]).
\end{equation}
It is easy to verify that
$\Omega^m:= \sum_{\p \in \G^m} \C \omega (\p) $ 
becomes an $\AS$-comodule via
\begin{equation}
\label{comOdef}
 \omega (\q) \mapsto 
 \sum_{\p \in \G^m} \omega (\p) \otimes e \binom{\p}{\q}  
 \quad (\q \in \G^m)
\end{equation}
for each $m \geq 0$.

For each $m \geq 0$, we set
\begin{equation}
 B \Omega^m =
 \bigl\{ (\lambda,\, 
 \lambda + \sum_{k \in I} \, \hat{i}\, ) 
 \in \V^2 \bigm|
 I \subset \{ 1, \ldots, N \},\,
 \mathrm{card} (I) = m \bigr\}.
\end{equation}  
Also we define 
$\Len\!: \G^m \to \Z_{\geq 0}$
by
\begin{equation}
 \EuScript{L}
 (\lambda\, |\, i_1, \ldots, i_m)
 =
 \mathrm{Card} \{ (k,l) |
 1 \leq k < l \leq N, i_k < i_l \} 
\end{equation}  
%
%
\begin{prop}
 For each $(\lambda,\mu) \in B \Omega^m$, 
 $\G^m_{\lambda \mu} \not= \emptyset$
 and $\omega_m (\lambda, \mu)$ $:=$ 
 $\sgn{\p} \omega (\p)$ does not 
 depend on the choice of 
 $\p \in \G^m_{\lambda \mu}$.
 Moreover 
 $\left\{ \omega_m (\lambda, \mu) \vert
 (\lambda, \mu) \in B \Omega^m \right\}$
 is a basis of $\Omega^m$. 
 In particular,  
 $\Omega^m = 0$ if $m > N$. 
 \end{prop}
\begin{pf}
 We will prove this lemma by means of Bergman's 
 diamond lemma \cite{Bergman}, 
 or rather its obvious generalization to
 the quotient algebras of $\C \langle \G \rangle$, 
 where $\langle \G \rangle = \coprod_m \G^m$.
 We define a ``reduction system'' 
 $S = S_1 \coprod S_2 \subset 
 \langle \G \rangle \times \C \langle \G \rangle$
 by setting
 \begin{gather*} 
 S_1 = \bigl\{ (\p, - \epsilon\, \p^{\dag}) \bigm| 
 \p = (\lambda\, |\, i,j) \in \G^2 [\searrow],\, i < j \bigr\},
 \\
 S_2 = \bigl\{ (\p, 0) \bigm| 
 \p = (\lambda\, |\, i_1, \ldots, i_m) \in \G^m,\,  m \geq 2,\, 
 \mathrm{card} \{ i_1, \ldots, i_m \} < m \bigr\}.
 \end{gather*}
 It is straightforward to verify that 
 the quotient
 $\C \langle \G \rangle / 
 \langle W - f \, \vert \, (W,f) \in S \rangle$
 is isomorphic to $\Omega$
 and that all ambiguities of $S$ are 
 resolvable. 
 Next, we introduce a semigroup partial order $\leq$ on 
 $\langle \G \rangle$ 
 by setting
 $(\lambda\, |\, i_1, \ldots, i_m) < 
 (\lambda\, |\, j_1, \ldots, j_n) $
 if either $m < n$, or $m = n$ and 
 $i_1 = j_1$, $\ldots$, $i_{k-1} = j_{k-1}$, 
 $i_k > j_k$ for some $1 \leq k \leq m$.
 Then $\leq$ is compatible with $S$ 
 and satisfies the descending chain condition.
 This completes the proof of the proposition.
\end{pf}
For each $0 \leq m \leq  N$,
we set 
$\Lambda_m = \hat{1}+ \cdots + \hat{m}$.
As an immediate consequence of 
\eqref{decompform} and the
the proposition above,			
we obtain the following result.
\begin{prop}
\label{O=L}
 For each $0 \leq m \leq  N$,
 we have $\Omega^m \cong L_{(\Lambda_m, m)}$
 as $\AS$-comodules. 
\end{prop}
Now we define the ``determinant'' 
${\det} = \sum_{\lambda, \mu \in \V} \det \binom{\lambda}{\mu}$ 
of $\AS$ to be the group-like element which corresponds to
the group-like comodule $\Omega^N$ and its basis
$\bigl\{ \bar{\omega} (\lambda)\bigr\}$ 
via Lemma \ref{glcom}, where					
$\bar{\omega} (\lambda) = 
D (\lambda) \omega_{N} (\lambda, \lambda)$ and
\begin{equation}
 D(\lambda) =
 \prod_{1 \leq i < j \leq N} 
 \frac{[d(\lambda\, |\, i,j)]}{[d(0\, |\, i,j)]}
 \quad (\lambda \in \V).
\end{equation}
Explicitly, we have
\begin{equation}
 \label{detformula}
 \det \binom{\lambda}{\mu} =
 \frac{D(\mu)}{D(\lambda)}
 \sum_{\p \in \G^N_{\lambda \lambda}}
 (- \epsilon)^{\Len (\p) + \Len (\q)}
 e \binom{\p}{\q},
\end{equation}
where $\q$ denotes an arbitrary 
element of $\G^N_{\mu \mu}$.
By \eqref{D(g)} and \eqref{gee} for $g = \det$, the quotient 
\begin{equation}
 \SS:= \AS / (\det - 1)
\end{equation}
naturally becomes a $\V$-face algebra,
which we call
$SU(N)_{L}$-{\it SOS algebra}. \\
The proof of the following lemma will be given in \S 8 and \S 13.
\begin{lem}
\label{cformula}
 For each $\p \in \G^1$, we have
 \begin{align}
  \label{c_OmegaNOmega1}
  c_{\Omega^N \Omega^1}
  \left( \bar{\omega} (\st (\p)) \otimes 
  \omega (\p) \right) =
  \epsilon^{N-1} \zeta^{- N} t
  \omega (\p) \otimes
  \bar{\omega} (\en (\p)), \\
  c_{\Omega^1 \Omega^N}
  \left( \omega (\p) \otimes 
  \bar{\omega} (\en (\p)) 
  \right) =
  \epsilon^{N-1} \zeta^{- N} t
  \bar{\omega} (\st (\p)) 
  \otimes \omega (\p). 
 \end{align}
\end{lem}
\begin{prop}
\label{SSisCQT}
 The element $\det$ belongs to the center
 of $\AS$. Moreover, if $\zeta$ satisfies
 \begin{equation}
  \label{etaeq}
  \zeta^N =
  \epsilon^{N-1} t,
 \end{equation}
 then 
 \begin{equation}
  \label{R(det-1,a)}
  \Rpm (\det -1, a) = 0 = \Rpm (a, \det -1)
  \quad (a \in \ASe).
 \end{equation}
 Hence, $\SS$ naturally becomes a 
 quotient CQT $\V$-face algebra of $\ASe$.  
\end{prop}
\begin{pf} (cf. \cite{gd})
 By computing the coaction of $\AS$ on 
 $\omega (\p) \otimes \bar{\omega} (\en (\p))$,
 in two ways via \eqref{c_OmegaNOmega1}, we obtain the first assertion. 
 We show the first equality of \eqref{R(det-1,a)}
 for $\pm = +$ and $a = e \binom{\p}{\q}$
 $(\p, \q \in \G^m, m \geq 0)$.
 By \eqref{Rpeeee} and \eqref{Rpeea}, it suffices to show
  \begin{equation}	
    \label{R(det,e)}
  \Rp \left( \det \binom{\en (\p)}{\st (\p)},\, e \binom{\p}{\q} \right)
  = \delta_{\p \q}
  \quad (\p, \q \in \G^m).
 \end{equation}
 For $m = 0$, this follows from \eqref{Rpeea} and \eqref{Rpeeee}.
 By computing the left-hand side of \eqref{c_OmegaNOmega1} via
 \eqref{cdef}, we obtain \eqref{R(det,e)} for $m = 1$.
 For $m \geq 2$, \eqref{R(det,e)} follows from 
 \eqref{Rp(a,bc)} and \eqref{D(g)} for $g = \det$
 by induction on $m$.
\end{pf}
				
Since the braiding of $\SS$ depends on the choice 
of the discrete parameter $\zeta$ satisfying \eqref{etaeq}, 
we sometimes write $\SSe$ instead of $\SS$. \\

To state our first main result, we recall the 
{\it fusion rule of} $SU(N)_L$-$WZW model$
in conformal field theory.
By \cite{GoodmanNakanishi}, it is characterized as 
the structure constant $N_{\lambda \mu}^{\nu}$
of a commutative $\Z$-algebra $\EuScript{F}$
(called the {\it fusion algebra of} $SU(N)_L$-$WZW model$)
with free basis $\{ \chi_{\lambda} \}_{\lambda \in \V}$
(i.e.,  $\chi_{\lambda} \chi_{\mu}$ $=$ 
$\sum_{\nu \in \V} N_{\lambda \mu}^{\nu} \chi_{\nu}$)
such that
\begin{equation}
\label{NlLmm} 
 N^{\mu}_{\lambda \Lambda_m} =
 \begin{cases} 
 1 & (\lambda, \mu) \in B \Omega^m \\ 
 0 & \text{otherwise} 
 \end{cases} 
\end{equation}
for each $\lambda, \mu \in \V$ and
$0 \leq m < N$.  
See Kac \cite{Kac} or Walton \cite{Walton} for an explicit
formula of $N_{\lambda \mu}^{\nu}$.

\begin{thm}[\cite{wzw}]
\label{SSrep}
 \rom{(1)} For each $\lambda \in \V$, up to isomorphism
 there exists a unique 
 right $\SS$-comodule $L_{\lambda}$
 such that 
 \begin{equation} 
 \label{dimLl0m}
  \dim L_{\lambda}(0, \mu) = \delta_{\lambda \mu}
  \quad (\mu \in \V).
 \end{equation} 
 Moreover, we have 
 \begin{equation} 
  \SS \cong 
  \bigoplus_{\lambda \in \V} \mathrm{End} (L_{\lambda})^*
 \end{equation} 
 as coalgebras.
 In particular, $L_{\lambda}$ is irreducible for each $\lambda \in \V$\\		

 \rom{(2)} The corepresentation ring 
 $K_0 ( \bold{Com}_{\SS} )$ is identified with the fusion algebra
 $\EuScript{F}$ of $SU(N)_L$-WZW model via 
 $\chi_{\lambda} = [ L_{\lambda} ]$.
 That is, we have
 \begin{equation} 
 \label{L0=1}
  [L_{0}] = 1,
 \end{equation} 
 \begin{equation} 
  [L_{\lambda}] [L_{\mu}] =
  \sum_{\nu \in \V} N^{\nu}_{\lambda \mu} [L_{\nu}].
 \end{equation} 
 \rom{(3)}
 We have
 \begin{equation} 
  \dim \left( L_{\nu} (\lambda, \mu) \right) = 
  N^{\mu}_{\lambda \nu}. 
 \end{equation} 
\end{thm} 
\begin{pf}							
 It is easy to verify that 
 $(\lambda, m) \in \boldsymbol\Lambda_{\G_{N,L}}$ if and only if 
 $m \in | \lambda | + N \Z_{\geq 0}$.
 Since $\C \V \det^n \bar{\otimes} L_{(\lambda,m)}$ $\cong$
 $L_{(\lambda,m + Nn)}$ by \eqref{MNlm}  			
 and \eqref{dimLlm*m}, we see that $\det$ satisfies 
 all conditions of Lemma \ref{sympredgle}, where 
 $\bar{\boldsymbol\Lambda} = \V$ and 
 $\varphi (\lambda, n) = (\lambda, | \lambda | + Nn)$.
 Therefore, we have Part (1), \eqref{L0=1} and 
 \begin{equation} 
  [L_{\lambda}] [L_{\mu}] =
  \sum_{\nu \in \V} N^{\nu}_{\lambda \mu} 
  (| \mu |)[L_{\nu}],
 \quad
  \dim \left( L_{\nu} (\lambda, \mu) \right) = 
  N^{\mu}_{\lambda \nu} (m)
  \quad ( (\nu, m) \in \boldsymbol\Lambda_{\G_{N,L}}). 
 \end{equation} 
 In particular, $\tilde{N}^{\nu}_{\lambda \mu}$ 
 $:=$ $N^{\nu}_{\lambda \mu} (m)$
 does not depend on the choice of $m$.
 Using Proposition \ref{O=L}, \eqref{decompform} and 
 \eqref{dimLl0m}, we obtain
 \begin{equation} 
  [L_{\lambda}] [L_{\Lambda_m}] =
  \sum_{\mu \in \V} 
  \dim ( L_{\lambda} \bar{\otimes} \Omega^m) (0, \mu)
  [L_{\mu}]
  =
  \sum_{\mu \in \V} \dim
  \left( \Omega^m (\lambda, \mu) \right)
  [L_{\mu}].
 \end{equation} 
 Thus the numbers $\tilde{N}^{\nu}_{\lambda \mu}$
 satisfies the condition of 
 $N^{\nu}_{\lambda \mu}$ stated above.
\end{pf} 
\begin{prop}
 The element $\det$ is not a zero-divisor of $\AS$.
 In particular, we have 
 $\det \binom{\lambda}{\mu} \not= 0$
 for each $\lambda, \mu \in \V$.
 Moreover, we have 
 \begin{equation} 
 \label{GLEAS}
  \mathrm{GLE} \left( \AS \right) =  
  \bigl\{ \sum_{\lambda, \mu \in \V} 
  \frac{c (\lambda)}{c (\mu)} \, \ema_{\lambda} e_{\mu} 
  {\det}^m  \bigm| m \in \Z_{\geq 0},\, 
  c(\lambda) \in \C^{\times} \,\, (\lambda \in \V) \bigr\}.
 \end{equation} 
\end{prop} 
\begin{pf}
 The first assertion follows from the fact that $\det$
 is simply reducible (see the proof of the theorem above).
 By Theorem \ref{Coststrthm} (1), every group-like comodule of $\AS$
 is isomorphic to $(\Omega^N)^{\overline{\otimes} m}$ for some $m \geq 0$.
 Hence the second assertion follows from Lemma \ref{glcom}.
\end{pf} 

\section{Module structure of $\Omega$}

%

Let $\HH$ be a CQT $\V$-face algebra over $\Bbb{K}$.
As in case $\HH$ is a CQT bialgebra, 
the correspondence $a \mapsto \Rp(\,,a)$
defines an antialgebra-coalgebra map from $\HH$
into the dual face algebra $\HH^{\circ}$ (cf. \cite{fa}).
Let $W$ be a right $\HH$-comodule.
Combining the above map with the left action 
\eqref{Xu} of $\HH^*$ on $W$, 
we obtain a right action of 
$\HH$ on $W$ given by 
\begin{equation}
\label{wa} 
 wa = \sum_{(w)} w_{(0)} \Rp (w_{(1)}, a)
 \quad (w \in W, a \in \HH).      
\end{equation}
Let $V$ be another $\HH$-comodule. Then we have
\begin{equation}
\label{(vw)a}
 \left( v \otimes w \right) a =
 \sum_{(a)} v a_{(1)} \otimes w a_{(2)} 
 \quad (v \in V(\lambda,\nu), w \in W(\nu,\mu), a \in \HH).
\end{equation}
If $\HH$ has an antipode and $W$ is finite-dimensional, 
then  we have 
\begin{equation}
\label{<va,w>}
 \langle v a, w \rangle =
 \langle v, w S^{-1} (a) \rangle
 \quad (v \in W \spcheck, w \in W),
\end{equation}
by \eqref{leftdualdef} and \eqref{Sid(R)}.

In this section, we give an explicit description of 
the right $\AS$-module structure of $\Omega$.
By \eqref{comOdef} and \eqref{UnivR=w}, 
we obtain the following.
\begin{lem}
For each $\s \in \G^m$ and
$\p, \q \in \G^n$ $(n \geq 0)$, we have
\begin{equation}
\label{osepq}
 \omega (\s) \, e \binom{\p}{\q} = 
 \sum_{\r \in \G^m}  
 w_{N,t,\epsilon} \!
 \left[ \p \frac[0pt]{\s}{\r} \q \right]\, 
 \omega (\r) 
 \quad (\s \in \G^m, \p, \q \in \G^n, n \geq 0).
\end{equation}
In particular, we have 
\begin{equation}
\label{Olmepq}
 \Omega (\lambda, \mu) e \binom{\p}{\q} \in  
 \delta_{\lambda, \st (\p)} \delta_{\mu, \st (\q)}\,
 \Omega (\en (\p), \en (\q)) .
\end{equation}
\end{lem}
\begin{lem}
\label{actOmega}				
 Let $(\lambda, \mu)$ be an element of $B \Omega^m$ $(m > 0)$
 and $\p$ $=$ $(\lambda\, |\, i_1, \ldots, i_m)$ an element of
 $\G^m_{\lambda \mu}$.
 Define the set $I$ and
 $C(\lambda\, |\, k,l) \in \C$ $(k \not= l)$ by
 $I = \{ i_1, \ldots, i_m \}$ and 
 \begin{equation}
  C(\lambda\, |\, k, l) = 
  \frac{[d(\lambda\, |\, k, l) + 1]}{[d(\lambda\, |\, k, l)]}
 \end{equation} 
 respectively. 
 Then for each $(\lambda\, |\, i), (\mu\, |\, j) \in \G^1$,
 we have\rom: 
 \begin{multline}
 \label{iinjnin}
  \omega (\p)
  e\! \binom{\lambda\, |\, i}{\mu\, |\, j}  
  \, = \,
  (- \zeta)^{- m} t^{-d(\lambda |\, i,j)} 
  \frac{1}{[d(\lambda |\, i,j)]} 
  \prod_{k \in I \setminus \{ i \} } 
  C(\lambda\, |\, i,k)\,
  \omega (\lambda + \hat{i}\, |\, i_2, \ldots, i_m, j)
  \\ 
  (i = i_1, j \not\in I), 
 \end{multline}
 \begin{multline}
 \label{iinjin}
  \omega (\p)
  e \binom{\lambda\, |\, i}{\mu\, |\, j}  
  \; = \;
  - \delta_{ij}
  (- \zeta)^{- m} t
  \prod_{k \in I \setminus \{ i \}} 
  C(\lambda\, |\, i,k) 
  \, \omega_m (\lambda + \hat{i}\, |\, i_2, \ldots, i_m, i) \\
 (i = i_1, j \in I), \qquad
 \end{multline}
 \begin{multline}
 \label{ininjnin}
  \omega_m (\p)
  e \binom{\lambda\, |\, i}{\mu\, |\, j}  
  \; = \;
  \delta_{ij}
  (\epsilon \zeta^{- 1})^m
  \prod_{k \in I} C(\lambda\, |\, i,k) 
  \, \omega_m (\lambda + \hat{i}\, |\, i_1, \ldots, i_m) \\ 
  (i, j \not\in I), \qquad\qquad
 \end{multline}
 \begin{equation}
 \label{ininjin}
  \qquad
  \omega_m (\p)    
  e \binom{\lambda\, |\, i}{\mu\, |\, j}  
  \; = \; 0
  \qquad\quad
  (i \not\in I, j \in I).
 \end{equation}
\end{lem}
\begin{pf}			
 These formulas are proved by induction on $m$ 
 in a similar manner.
 Here we give the proof of \eqref{iinjnin}.	
 Since $\Omega^m$ is a quotient module of 
 $\C \G^{m-1} \bar{\otimes} \C \G^1$, 		
 the left-hand side of \eqref{iinjnin} is rewritten as
 $\sum_q A_q B_q$ with
 \begin{equation}
  A_q = \omega_{m-1} (\lambda, \nu)    
  e\! \binom{\lambda\, |\, i}{\nu\, |\, q},
 \quad
  B_q = \omega_1 (\nu, \mu)    
  e\! \binom{\nu\, |\, q}{\mu\, |\, j}
 \end{equation}
 by \eqref{(vw)a}, where 
 $\nu = \lambda + \hat{i}_1 + \cdots + \hat{i}_{m-1}$
 and the summation is taken over for all $1 \leq q \leq N$
 such that $(\nu\, |\, q) \in \G^1$.			
 Since $(\nu + \hat{q}, \mu + \hat{j}) \in \G^1$
 only if $q = i_m$ or $j$, we have
 \begin{equation}
 \label{oe=AB+AB}
 \omega_m (\lambda, \mu) 
 e\! \binom{\lambda\, |\,i}{\mu\, |\,j} 
 =
 \begin{cases}
  A_{i_m} B_{i_m} + A_j B_j 
  & (\nu\, |\, j) \in \G^1\\
  A_{i_m} B_{i_m} 
  & \mathrm{otherwise.}
 \end{cases}
 \end{equation}
 Using inductive assumption, we see that
 the right-hand side of \eqref{oe=AB+AB} equals
 \begin{multline}
 \label{iinjninpre} 
  (- \zeta)^{- m} t^{-d(i,j)} 
  \left(
  \frac{1}{[d(i, i_m)] [d(i_m, j)]} +
  \frac{[d(i_m, j) - 1]}{[d(i, j)] [d(i_m, j)]} 
  \right) \\
  \prod_{n=2}^{m-1} 
  C(\lambda\, |\, i, i_n)\,
  \omega_m (\lambda + \hat{i}, \mu + \hat{j})
  \qquad\qquad
 \end{multline}
 if $(\nu\, |\, j) \in \G^1$,
 where $d(k, l) = d (\lambda\, |\, k, l)$.
 Applying $[a + b + 1] + [a][b] =$
 $[a + 1][b + 1]$ to $a = d(i, i_m)$
 and $b = d(i_m, j) -1$, we see that 
 \eqref{iinjninpre} equals the right-hand side
 of \eqref{iinjnin}
 (even if $(\nu\, |\, j) \not\in \G^1$).
 Next suppose that $(\nu\, |\, j) \not\in \G^1$.
 It suffices to verify that the second term
 in the parentheses of \eqref{iinjninpre} is zero.
 In case $j = 1$, we obtain $\nu_1 = L$.
 Using this together with $1 \not\in I$, we see that
 $\lambda_1 = L$.
 On the other hand, since $(\mu\, |\, 1) \in \G^1$, 
 we have $L - 1 \geq \mu_1 = L - \delta_{i_m N}$.
 Hence, $i_m = N$ and  
 $d(i_m, j) -1 = - L - N$.
 In case $j > 1$, we obtain
 $d(i_m, j) -1 = 0$ in a similar manner.
 Thus we complete the proof of \eqref{iinjnin}. 
\end{pf}			
The following lemma is frequently used in the sequel.
\begin{lem}
\label{irrG1}
 As a right $\AS$-module, $\C \G^1 \cong \Omega^1$
 is irreducible. Hence $\C \G^1$ is also irreducible 
 as a left $\AS^*$-module. 
\end{lem}
\begin{pf}				
 Let $W$ be a non-zero submodule of $\C \G^1$.
 Since $W = \bigoplus_{\lambda \mu} W \ema_{\lambda} e_{\mu}$
 and
 $W \ema_{\lambda} e_{\mu} \subset \C \G^1_{\lambda \mu}$, 
 we have $\s_0 \in W$ for some 
 $\s_0 = (\lambda\, |\, i) \in \G^1$.
 To show $\C \G^1 = \s_0\, \AS$, we introduce the oriented graph 
 $\EuScript{H}$ determined by $\EuScript{H}^0 = \G^1$ and
 \begin{equation}
 \mathrm{card} \left( \EuScript{H}^1_{\s \r} \right) =
 \begin{cases}
  1 & \s\, e\! \binom{\p}{\q} \in \C^{\times} \r 
  \quad (\exists \p, \q \in \G^1) \\
  0 & \mathrm{otherwise}.
 \end{cases}
 \end{equation}
 It suffices to show that 		
 $\langle \EuScript{H} \rangle_{\s_0\, \r} \not=$
 $\emptyset$ for every $\r \in \G^1$, where
 $\langle \EuScript{H} \rangle_{\s\, \r} =$
 $\cup_m \EuScript{H}^m_{\s\, \r}$.
 We note that 
 \begin{equation}
 \label{H1necond3}
  \EuScript{H}^1_{(\mu\, |\, j)\, (\mu + \hat{j}\, |\, k)} 
  \not= \emptyset
  \quad \mathrm{if}\;\,	
  (\mu\, |\, j), (\mu + \hat{j}\, |\, k)
  \in \G^1\; \mathrm{and}\; j \not= k,
 \end{equation}
 \begin{equation}
 \label{H1necond1}			
  \EuScript{H}^1_{(\mu\, |\, j)\, (\mu + \hat{k}\, |\, j)} 
  \not= \emptyset
  \quad \mathrm{if} \;\, 			
  (\mu\, |\, j), (\mu + \hat{k}\, |\, j)
  \in \G^1\; \mathrm{and}\; j \not= k,
 \end{equation}
 \begin{equation}
 \label{H1necond2}
  \EuScript{H}^1_{(\mu\, |\, j)\, (\mu + \hat{j}\, |\, j)} 
  \not= \emptyset
  \quad \mathrm{if}\;\, 
  (\mu\, |\, j, j) \in \G^2.
 \end{equation} 
 Using \eqref{H1necond1}, we obtain 
 $\langle \EuScript{H} \rangle_{\s_0\, \s_1}$
 $\not=$ $\emptyset$, 
 where $\s_1 =$
 $(L \Lambda_{i-1} +			
 \sum_{k=i}^{N-1} \lambda_k \hat{k}\,|\, i)$.
 Suppose $i \not= N$. 
 Using \eqref{H1necond2}, and then 
 using \eqref{H1necond1}, we obtain 
 $\langle \EuScript{H} \rangle_{\s_1\, \s_2}$
 $\not=$ $\emptyset$, where 	
 $\s_2 =$
 $((L-1) \Lambda_{N-1} + \Lambda_{i-1}\, |\,i)$.
 Using \eqref{H1necond3} and \eqref{H1necond1} 
 respectively, we obtain
 $\langle \EuScript{H} \rangle_{\s_2\, \s_3},$
 $\langle \EuScript{H} \rangle_{\s_4\, (0 | 1)}$
 $\not=$ $\emptyset$ and 
 $\langle \EuScript{H} \rangle_{\s_3\, \s_4}$
 $\not=$ $\emptyset$ respectively, where
 $\s_3 =$
 $((L-1) \Lambda_{N-1} + \Lambda_{N-2}\, |\,N-1)$.
 and $\s_4 =$ $(\Lambda_{N-2}\, |\, N-1)$.
 Therefore we obtain 
 $\langle \EuScript{H} \rangle_{\s_0 (0|1)}$
 $\not=$ $\emptyset$ if $i \not= N$.
 By similar consideration, we also obtain
 $\langle \EuScript{H} \rangle_{\s_0 (0|1)}$
 $\not=$ $\emptyset$ in case $i = N$, 
 and also,
 $\langle \EuScript{H} \rangle_{(0|1) \r}$
 $\not=$ $\emptyset$ for every $\r \in \G^1$.
 Thus, we have verified the first assertion.
 The second assertion is obvious since the image of
 $a \mapsto \Rp (\, , a)$ is a subalgebra of $\AS^*$.
\end{pf}
Now we begin to prove Lemma
\ref{cformula}.
By \eqref{Olmepq} and Lemma \ref{actOmega}, we have
\begin{equation}
 \bar{\omega} (\nu)    
 e \binom{\lambda\, |\, i}{\mu\, |\, j}  
 \; = \;
 \epsilon^{N-1} \zeta^{- N} t
 \delta_{ij} \delta_{\nu \lambda} \delta_{\nu \mu}
 \bar{\omega} (\nu + \hat{i}).    
\end{equation}
Using \eqref{(vw)a} and this equality, we see that both 
$\Omega^N \bar{\otimes} \C \G^1 \to \C \G^1$;
$\bar{\omega} (\st (\p)) \otimes \p$
$\mapsto$
$\p \otimes \bar{\omega} (\en (\p))$
and 				
$\C \G^1 \bar{\otimes} \Omega^N \to \C \G^1$;
$\p \otimes \bar{\omega} (\st (\p))$
$\mapsto$
$\p \otimes \bar{\omega} (\en (\p))$
are isomorphisms of right $\AS$-modules.
Hence, by Lemma \ref{irrG1} and 
Schur's Lemma, we have
\begin{equation}
\label{cONO1=theta}
 c_{\Omega^N \Omega^1}  
 (\bar{\omega} (\st (\p)) \otimes \p) 
 = \vartheta\,
 \p \otimes \bar{\omega} (\en (\p))
 \quad (\p \in \G^1)
\end{equation}
for some constant $\vartheta$. We will prove 
$\vartheta = \epsilon^{N-1} \zeta^{- N} t$ in \S12.

\section{transposes and complex conjugates}
The following proposition is an immediate consequence 
of the following reflection symmetry:
\begin{equation}
 w_{N,t,\epsilon} \!
 \left[ \r \frac[0pt]{\p}{\s} \q \right] =
 \left(
 \frac{\kappa (\r \cdot \s)}{\kappa (\p \cdot \q)} 
 \right)^2
 w_{N,t,\epsilon} \!
 \left[ \p \frac[0pt]{\r}{\q} \s \right], 
\end{equation} 
where $\kappa$ is as in \eqref{kappadef}. 
\begin{prop}
\label{trA}		
 There exists an algebra-anticoalgebra
 map $\AS \to \AS;$ $a \mapsto a^{\tenc}$
 given by
 \begin{equation}
 \label{tencdef}
  e \binom{\p}{\q}^{\tenc} = 
  \left( \frac{\kappa (\p)}{\kappa (\q)} \right)^2
  e \binom{\q}{\p}
  \quad (\p, \q \in \G^m, m \geq 0).
 \end{equation}
 Moreover it satisfies $(a^{\tenc})^{\tenc}$
 $=$ $a$ and
 \begin{equation}
 \label{R&tenc}
  \Rpm  \left(a^{\tenc},\, b^{\tenc} \right) =
  \Rpm  \left(b,\, a \right)
 \end{equation}
 for each  $a, b \in \AS$ and $\zeta \in \C^{\times}$.
\end{prop}
The following proposition is needed to construct the ``cofactor 
matrix.''
\begin{prop}
\label{trS}
 The element $\det$ satisfies
 ${\det}^{\tenc} = \det$.
 Hence, $\tenc$ induces an algebra-anticoalgebra 
 involution of $\SS$, which satisfies \eqref{R&tenc}.
\end{prop}
\begin{pf}
 Since ${\det}^{\tenc}$ is a group-like element 
 of $\AS$, we have 
 \begin{equation}
 {\det}^{\tenc} = g \det;
 \quad
 g = \sum_{\lambda, \mu \in \V} 
 \frac{c(\lambda)}{c(\mu)} \, \ema_{\lambda} e_{\mu}
 \end{equation}
 by \eqref{GLEAS}, where $c(\lambda)$ $(\lambda \in \V)$ denotes
 some nonzero constant.
 Since both $\det$ and ${\det}^{\tenc}$ are central 
 and $\det$ is not a zero divisor, $g$ is central.
 Hence, by Lemma \ref{irrG1} and Schur's lemma, 
 we have $\p \, g = c \, \p$ $(\p \in \G^1)$ for some $c \in \C$.
 Hence we have $c(\lambda) = c^{| \lambda |} c(0)$ 
 for each $\lambda \in \V$.
 In order to prove $c = 1$, 
 we compute ${\det \binom{0}{\hat{1}}}^{\tenc}$
 in two ways. 
 Using
 \begin{equation}
  [d(\p) - 1] \, e \! \binom{\p}{\q} = 
  - \epsilon [d(\p) + 1] \, e \! \binom{\p^{\dag}}{\q}
  \quad (\p \in \ {\G}^2[\;\searrow\;], \,
  \q \in  {\G}^2[\;\downarrow\;]),  
 \end{equation}
 we obtain
 \begin{equation}
  \sum_{i=1}^{k} \sgn{\p_i} 
  e \binom{\p_i}{0\, |\, 1, \ldots, N} =
  \sgn{\p_k} [k]^2 \, 
  e \binom{\p_k}{0\, |\, 1, \ldots, N} 
 \end{equation}
 by induction on $k$, where
 $\p_i = (\hat{1}\, |\, 2, \ldots, i, 1, i+1, \ldots, N)$.
 Substituting $k = N$ in this equality, we get
 \begin{align*}
  \det {\binom{0}{\hat{1}}}^{\tenc} & = \,\,
  c^{|\hat{1}| - |0|}
  \det \binom{\hat{1}}{0} 
  \\ 
  & =
  (- \epsilon)^{N(N-1)/2 + \Len{(\p_N)}}
  c [N]\,  e \! 
  \begin{pmatrix}
   \hat{1}\,|\!\!\! & 2,\!\!\! & 3,\!\!\! & \ldots,\!\!\! & N,\!\!\!     &  1 \\
    0\,|\!\!\!      & 1,\!\!\! & 2,\!\!\! & \ldots,\!\!\! & N - 1,\!\!\! & N \\
  \end{pmatrix}.
 \end{align*}
 On the other hand, using \eqref{detformula}
 and \eqref{tencdef}, we see that the right-hand side of the 
 above equality agrees with $c \,{\det \binom{0}{\hat{1}}}^{\tenc}$.
 This completes the proof of the proposition.
\end{pf}
Next suppose 
$t = \exp (\pm \frac{\pi i}{N + L})$, or
$- \exp (\pm \frac{\pi i}{N + L})$ with $N + L \in 2 \Z$.
Then, we have $\kappa (\p) > 0$ for each 
$\p \in \G^m$ $(m > 0)$.
Moreover, for each $\zeta$ with $| \zeta | = 1$, 
the Boltzmann weight $w_{N,t}$ satisfies
\begin{equation}
\label{wNt*}
 w_{N,t,\epsilon}^{-1} \!
 \left[ \r \frac[0pt]{\p}{\s} \q \right] =
 \left(
 \frac{\kappa (\r \cdot \s)}{\kappa (\p \cdot \q)} 
 \right)^2 \,
 \overline{w_{N,t,\epsilon} \!
 \left[ \p \frac[0pt]{\r}{\q} \s \right]}.
\end{equation} 
Similarly to Proposition \ref{trA} and 
Proposition \ref{trS}, we obtain the following.
\begin{prop}				
\label{SSiscpt}
 Set $t = \pm \exp (\pm \frac{\pi i}{N + L})$
 if $N + L \in 2 \Z$, and
 $t = \exp (\pm \frac{\pi i}{N + L})$ if
 $N + L \in 1 + 2 \Z$. 
 Then for each solution $\zeta$ of \eqref{etaeq},
 both $\SSe$ and $\ASe$ are compact CQT $\V$-face algebras 
 of unitary type with costar structure 
 \begin{equation}
 \label{costardef}
  e \binom{\p}{\q}^{\times} = 
  \left( \frac{\kappa (\p)}{\kappa (\q)} \right)^2
  e \binom{\q}{\p}
  \quad (\p, \q \in \G^m, m \geq 0).
 \end{equation}
\end{prop}
In fact, both $\SS$ and $\AS$ are spanned by unitary
matrix corepresentations				
$[e_u \binom{\p}{\q}]_{\p, \q \in \G^m}$ 
$(m \geq 0)$ given by
 \begin{equation}
 \label{esymdef}
  e_u \binom{\p}{\q} = 
  \frac{\kappa (\q)}{\kappa (\p)} \,
  e \binom{\p}{\q}
  \quad (\p, \q \in \G^m, m \geq 0).
 \end{equation}

\section{Antipodes and ribbon functionals}
\begin{lem}
The $\V$-face algebra $\SS$ has an antipode given by
\begin{equation}
\label{Santipode}
 S \left( e \binom{\lambda\,|\,i}{\mu\,|\,j} \right) = 
 (- \epsilon)^{i+j} 
 \frac{D(\lambda + \hat{i})}{D(\mu + \hat{j})}
 \sum_{\aaa} 
 (- \epsilon)^{\Len (\aaa) + \Len (\bbb)} e \binom{\aaa}{\bbb},
\end{equation}
%
%
where $\bbb$ denotes an arbitrary element of 
$\G^{N-1}_{\lambda + \hat{i}\,\lambda}$
and the summation is taken over all 
$\aaa \in \G^{N-1}_{\mu + \hat{j}\, \mu}$.
Moreover, we have 
\begin{equation}
\label{S2formula}
  S^2 \left( e \binom{\p}{\q} \right) = 
 \frac{D(\en (\p))D(\st(\q))}{D(\st(\p))D(\en(\q))} \, 
 e \binom{\p}{\q}
 \quad (\p, \q \in \G^m, m \geq 0).
 \end{equation}
\end{lem}
\begin{pf}
 (cf. \cite{TakeuchiMat}, \cite{qcg})
 Let $Y_{(\lambda\, |\,i)\,(\mu\, |\,j)}$ denote
 the right-hand side of \eqref{Santipode}
 viewed as an element of $\AS$.
 By \cite[\S 7]{cpt}, it suffices to verify that
 \begin{equation}
 \label{cofactor}
  \sum_{\r \in \G^1} Y_{\p \r} X_{\r \q} =   
  \delta_{\p \q} e_{\en (\p)} \det,
 \quad  
  \sum_{\r \in \G^1} X_{\p \r} Y_{\r \q} =   
  \delta_{\p \q} \ema_{\st (\p)} \det,  
 \end{equation}
 where
 $X_{\p \q}$ $=$ $e \binom{\p}{\q}$
 $(\p, \q \in \G^1)$.
 We define a basis 
 $\bigl\{ \bar{\omega} (\tilde{\p})
 \bigm| \p \in \G^1 \bigr\}$  
 of $\Omega^{N-1}$ by
 $\bar{\omega} (\lambda + \hat{i}, \lambda) =
 (- \epsilon)^{i-1} D(\lambda + \hat{i})
 \omega_{N-1} (\lambda + \hat{i}, \lambda)$, 
 where
 $\tilde{\p} = (\mu, \lambda) \in B \Omega^{N-1}$ for
 $\p = (\lambda, \mu) \in \G^1$.
 Then, the coaction of $\AS$ on $\Omega^{N-1}$ is given by
 $\bar{\omega} (\tilde{\q}) \mapsto
 \sum_{\p \in \G^1} \bar{\omega} (\tilde{\p}) \otimes Y_{\q \p}$
 and the multiplication of $\Omega$ gives maps
 \begin{gather*}
  \Omega^{N-1} \bar{\otimes} \Omega^1 \to \Omega^N;
  \quad 
  \bar{\omega} (\tilde{\q}) \otimes \omega(\p) \mapsto
  \delta_{\p \q} \bar{\omega} (\en (\p)), \\
  \Omega^1 \bar{\otimes} \Omega^{N-1} \to \Omega^N;
  \quad
  \omega(\p) \otimes \bar{\omega} (\tilde{\q}) \mapsto
  \delta_{\p \q} (- \epsilon)^{N-1} 
  \frac{D(\en (\p))}{D(\st (\p))}
  \bar{\omega} (\st (\p)). 
 \end{gather*}
 Since these maps are compatible with the coaction of $\AS$, 
 we have
 the first formula of \eqref{cofactor}
 and 
 \begin{equation}
  \sum_{\r} W_{\r \p} Y_{\q \r} =
  \delta_{\p \q} e_{\st (\p)} \det, 
 \end{equation}
 where $W_{\p \q} \in \AS$ denotes the right-hand side of 
 \eqref{S2formula}.
 Applying $\tenc$ to this equality, we obtain
 \begin{equation}
  \sum_{\r} X_{\p \r} Z_{\r \q} =
  \delta_{\p \q} \ema_{\st (\p)} \det;
  \quad 
  Z_{\p \q} =
  \frac{D(\en (\p)) D(\st (\q)) \kappa (\p)^2}
  {D(\st (\p)) D(\en (\q)) \kappa (\q)^2} \, 
  Y_{\q \p}^{\tenc}.
 \end{equation}
 Computing 
 $\sum_{\r \s} Y_{\p \r} X_{\r \s} Z_{\s \q}$
 in two ways, we obtain 
 $Y_{\p \q} =  Z_{\p \q}$.
 This proves the second equality of \eqref{cofactor}.
 Finally, Computing 
 $\sum_{\r \s} W_{\r \q} S ( X_{\s \r} ) S^2 (X_{\p \s})$
 in the algebra $\SS$
 in two ways, we obtain \eqref{S2formula}.
\end{pf}
\begin{prop}
\label{SSrib}
 For each $t$ and $\zeta$, $\SSe$ becomes a coribbon Hopf $\V$-face algebra,
 whose braiding, antipode $S$ and modified ribbon functional
 $\cal{M} = \cal{M}_{1}$ are given by \eqref{UnivR=w}, 
 \eqref{Santipode}
 and the following formulas\rom: 
 \begin{equation}
 \label{Mformula}
  \cal{M} \left( e \binom{\p}{\q} \right) =
  \delta_{\p \q} \,
  \frac{D(\en (\p))}{D(\st (\p))}
  \quad (\p, \q \in \G^m, m \geq 0).
 \end{equation}
 Moreover, we have 
 \begin{equation}
 \label{d=D}
  \dim_q (L_{\lambda}) = D(\lambda)
 \end{equation}
 for each $\lambda \in \V$.
 When $N$ is even, there exists another ribbon
 functional $\cal{M}_{-1}$ given by  
 \begin{equation}
 \label{M-formula}
  \cal{M}_{-1} \left( e \binom{\p}{\q} \right) =
  \delta_{\p \q} \,
  (-1)^{ | \en (\p) | - | \st (\p) | }
  \frac{D(\en (\p))}{D(\st (\p))}
  \quad (\p, \q \in \G^m, m \geq 0).
 \end{equation}
 The quantum dimension of the corresponding ribbon category
 is given by
 \begin{equation}
 \label{d=-D}
  \dim_q^{-1} (L_{\lambda}) = (-1)^{| \lambda |} D(\lambda).
 \end{equation}
 \end{prop} 

\begin{pf}
 Let $\cal{M} \in \HH^*$ be as in \eqref{Mformula}.
 Using \eqref{S2formula}, we see that 
 $\cal{M}$ satisfies \eqref{MXM-}.
 Hence, it suffices to verify that
 \begin{equation}
 \label{U1U2=M2}	
  (\cal{U}_1 \cal{U}_2) \left( e \binom{\p}{\q} \right)
  = \cal{M}^2 \left( e \binom{\p}{\q} \right)
  \quad (\p, \q \in \G^m).
 \end{equation}
 for each $m \geq 0$.
 By \eqref{UXU-} and \eqref{MXM-}, $\cal{U}_1 \cal{M}^{-1}$
 is a central element of $\AS^*$.
 Hence by Lemma \ref{irrG1} and Schur's lemma, 
 we have $\cal{U}_1 \p = \vartheta \cal{M} \p$
 $(\p \in \G^1)$ for some constant $\vartheta$.
 Using \eqref{Um}, \eqref{Sid(R)} and \eqref{S2formula}, 
 we compute
  \begin{align}
  \label{Uformula}		
   \cal{U}_1^{-1} \left( e \binom{0\, |\, 1}{0\, |\, 1}
   \right) = &
   \sum_{\r \in \G^1}			
   \Rp \left( S^2 \left( e \binom{\r}{0\, |\, 1} \right),\, 
   e \binom{0\, |\, 1}{\r} \right) \\
   = &
   \sum_{\nu = 2 \hat{1}, \hat{1} + \hat{2}}		
   \frac{D(0)D(\nu)}{D(\hat{1})^2}
   w_{N,t}
   \begin{bmatrix}
    0       & \hat{1} \\
    \hat{1} & \nu
   \end{bmatrix} \\
   = & \;
   \zeta^{-1} t^{N} \frac{1}{[N]}.
  \end{align}
 This shows that $\vartheta = \zeta t^{-N}$, 
 and similarly, we obtain 
 $\cal{U}_2 \p = \vartheta^{-1} \cal{M} \p$
 $(\p \in \G^1)$. This proves \eqref{U1U2=M2}  
 for $m = 1$.
 For $m \geq 2$, \eqref{U1U2=M2} follows from 
 \eqref{U1U2=M2} for $m = 1$ by induction on $m$,
 using the fact that both $\cal{M}^2$ and 
 $\cal{U}_1 \cal{U}_2$ are group-like 
 (cf. \eqref{U(ee)}, \eqref{m*(U)}).		
 The second assertion follows from \eqref{Mformula}
 and Lemma \ref{Trq}.
\end{pf}
We denote by $\SSi$ the coribbon Hopf $\V$-face algebra
$(\SSz, \cal{M}_{\iota})$,
where $\iota = \pm 1$ if $N$ is even and  
$\iota = 1$ if $N$ is odd.
\begin{lem}  
 If $t = \exp (\pm \frac{\pi i}{N + L})$, 
 or  $t = - \exp (\pm \frac{\pi i}{N + L})$
 and $N$ is odd, then we have $D (\lambda ) > 0$
 for every $\lambda \in \V$.
 If $t = - \exp (\pm \frac{\pi i}{N + L})$ 
 and $N$ is even, then we have 
 $(-1)^{| \lambda |} D (\lambda ) > 0$.
\end{lem} 
\begin{pf}  
 Straightforward.
\end{pf} 
\begin{prop}  
 When $t = \exp (\pm \frac{\pi i}{N + L})$, 
 or  $t = - \exp (\pm \frac{\pi i}{N + L})$
 and $N, L \in 1 + 2 \Z$, 
 the Woronowicz functional $\cal{Q}$ of 
 $\SS$ is given by \eqref{Mformula}.
 While when $t = - \exp (\pm \frac{\pi i}{N + L})$ 
 and $N, L \in 2 \Z$ , $\cal{Q}$ 
 is given by \eqref{M-formula}. 
\end{prop} 
\begin{pf}
 We will prove the first assertion. We set 
 \begin{equation}
  M = 			
  \left[ \cal{M} \left( 
  e_u \binom{\p}{\q}
  \right) \right]_{\p, \q \in \G^1}
 \quad 
  Q =			
  \left[ \cal{Q} \left( 
  e_u \binom{\p}{\q} \right) \right]_{\p, \q \in \G^1}
 \end{equation}
 where $\cal{Q}$ denotes the Woronowicz functional of $\SS$
 and $e_u \binom{\p}{\q}$ is as in 
 \eqref{esymdef}.					
 By \eqref{MXM-}, \eqref{QXQ-} and Lemma \ref{irrG1}, we have
 $M = \vartheta Q$ for some $\vartheta$.
 Since $M$ is positive by \eqref{Mformula} and 
 the lemma above, 
 we have $\vartheta > 0$.
 Since the quantum dimension satisfies 
 $\dim_q L = \dim_q L\spcheck$ for every $L$, we have
 $\mathrm{Tr} (M) = \mathrm{Tr} (M^{-1})$.
 By \eqref{TrQ}, this proves $M = Q$. Now the assertion
 $\cal{M} \bigl( e_u \binom{\p}{\q} \bigr)$ $=$
 $\cal{Q} \bigl( e_u \binom{\p}{\q} \bigr)$
 $(\p, \q \in \G^m)$ easily follows from the fact that both
 $\cal{M}$ and  $\cal{Q}$ are group-like, by induction on $m$.
\end{pf}		
%


\section{The modular tensor category}
Let $\cal{C}$ be a ribbon category,
which is additive over a field $\Bbb{K}$.
We say that $\cal{C}$ is {\it semisimple}
if there exist a set $\V_{\cal{C}}$, an involution 
${}^{\lor}\!: \V_{\cal{C}} \to \V_{\cal{C}}$, 
an element $0 \in \V_{\cal{C}}$
and simple objects $L_{\lambda}^{\cal{C}}$ 
$(\lambda \in \V_{\cal{C}})$
such that every object of $\cal{C}$ is isomorphic to
a finite direct sum of $L_{\lambda}^{\cal{C}}$'s, 
and that			
\begin{equation}
 L_0^{\cal{C}} \cong \bold{1},
\quad
 (L_{\lambda}^{\cal{C}})^{\lor}
 \cong
 L_{\lambda^{\lor}}^{\cal{C}},
\end{equation}
\begin{equation}
 \cal{C} (L_{\lambda}^{\cal{C}}, L_{\mu}^{\cal{C}}) =
 \begin{cases}
  \Bbb{K} & (\lambda = \mu) \\
  0 & (\lambda \ne \mu) \\
 \end{cases}
\end{equation}
for each $\lambda, \mu \in \V^{\cal{C}}$,
where $\bold{1}$ denotes the unit object of $\cal{C}$.
For a semisimple ribbon category $\cal{C}$, we define its
{\it fusion rule} $N^{\nu}_{\lambda \mu}$
$(\lambda, \mu, \nu \in \V_{\cal{C}})$ and 
{\it S-matrix} $S^{\cal{C}}$ $=$ 
$[S^{\cal{C}}_{\lambda \mu}]_{\lambda, \mu \in \V_{\cal{C}}}$
by 
\begin{gather}
 [L_{\lambda}^{\cal{C}}][L_{\mu}^{\cal{C}}] =  
 \sum_{\nu \in \V_{\cal{C}}} 
 N^{\nu}_{\lambda, \mu} [L_{\nu}^{\cal{C}}], \\
 S^{\cal{C}}_{\lambda \mu} 
 =
 \mathrm{Tr}_q (c_{L_{\mu}^{\cal{C}} L_{\lambda}^{\cal{C}}}
 \circ 
 c_{L_{\lambda}^{\cal{C}} L_{\mu}^{\cal{C}}}).
\end{gather}
By definition, we have 
\begin{equation}
\label{Sl0C}
 S^{\cal{C}}_{\lambda 0} =
 \dim_q L_{\lambda}^{\cal{C}}.
\end{equation}
Since the twist $\theta$ satisfies
\begin{equation}
\label{thetaVW}
 c_{W V} \circ c_{V W} =
 \theta_{V \otimes W} \circ 
 (\theta_V \otimes \theta_W)^{-1},
\end{equation} 
$S^{\cal{C}}$ satisfies
\begin{equation}
\label{SlmC}
 S_{\lambda \mu}^{\cal{C}} =
 \sum_{\nu \in \V_{\cal{C}}}
 \frac{\theta_{\nu}}{\theta_{\lambda} \theta_{\mu}} 
 N_{\lambda \mu}^{\nu} 
 \dim_q (L_{\nu}^{\cal{C}}),
\end{equation}
where $\theta_{\lambda} \in \Bbb{K}$ is defined by
$\theta_{L_{\lambda}^{\cal{C}}} =  
\theta_{\lambda}   
\mathrm{id}_{L_{\lambda}^{\cal{C}}}$.
Moreover, $S = S^{\cal{C}}$ satisfies the following 
{\it Verlinde's formula}
(cf. \cite{Verlinde,MooreSeiberg,Turaev}):
\begin{equation}
\label{Verlinde}
 {S_{\nu,0}}
 \sum_{\xi \in \V} N_{\lambda \mu}^{\xi} 
 S_{\xi \nu} 
 = S_{\lambda \nu} S_{\mu \nu},
\end{equation}
where $\lambda$, $\mu$ and $\nu$ denote arbitrary
elements of $\V = \V_{\cal{C}}$.
Let $\cal{C}$ be a semisimple ribbon category.
We say that $\cal{C}$ is a {\it modular tensor category}
(or {\it MTC}) if $\V_{\cal{C}}$ is finite and 
the matrix $S^{\cal{C}}$ is invertible.
If, in addition, $\cal{C}$ is unitary as a ribbon category, 
then it is called a {\it unitary MTC}.
It is known that each (unitary) MTC gives rise to a 
(unitary) 3-dimensional topological quantum field theory (TQFT), 
hence, in particular, an invariant of 3-manifolds
of Witten-Reshetikhin-Turaev type (cf. V. Turaev \cite{Turaev}). 

The most well-known example of MTC is obtained
as a certain semisimple quotient
$\cal{C} (\frak{g}, \kappa)$ of a category
of representations of the quantized enveloping 
algebra $U_q (\frak{g})$ of finite type
in the case when $q$ is a root of unity
\cite{Andersen,GelfandKazhdan,Kirillov,TuraevWenzl}.
When $\frak{g} = \frak{sl}_N$, the simple objects
$L_{\lambda}^U$ of $\cal{C} (\frak{sl}_N, N + L)$
$(L \geq 1)$ are also indexed by the set $\V = \V_{N L}$ 
given by \eqref{Vdef} and the fusion rules agree
with those of $SU(N)_L$-WZW model.
The quantum dimension and the constant $\theta_{\lambda}$
for $\cal{C} (\frak{sl}_N, N + L)$ 
are given by 

\begin{equation}
\label{dimqU}		
 \dim_q ( L_{\lambda}^U ) = D (\lambda)_{t_0},
\quad
 \theta_{\lambda} = \zeta_0^{\, (\lambda | \lambda + 2 \rho )^{\sim}}
\end{equation}
respectively, where 
$\zeta_0 = \exp (\frac{\pi i}{N (N + L)})$,
$t_0 = \zeta_0^N$, 
$\rho = \Lambda_1 + \cdots + \Lambda_{N-1}$, 
$(\,|\,)^{\sim} = N (\, |\, )$ 
and $(\,|\,)$ denotes the usual 
inner product of $\R^N$.
Moreover, the S-matrix of
$\cal{C} (\frak{sl}_N, N + L)$ is given by
$S^U =  S^{1} (\zeta_0)$.
Here, for each primitive $2 N (N + L)$-th root $\zeta$
of unity, we define the matrix 
$S^{\iota} (\zeta)$ by the following 
Kac-Peterson formula (cf. \cite{Kirillov}):
\begin{equation}
\label{SlmU} 
 S^{\iota} (\zeta)_{\lambda \mu} = 
 \iota^{| \lambda | + |\mu | }
 \frac{\sum_{w \in \frak{S}_N} (-1)^{l (w)} 
 \zeta^{-2 (w (\lambda + \rho)\, |\, \mu + \rho )^{\sim}}}
 {\sum_{w \in \frak{S}_N} (-1)^{l (w)} 
 \zeta^{-2 (w (\rho)\, |\, \rho)^{\sim}}},
\end{equation}
where $\iota = \pm 1$ if $N \in 2 \Z$,
$\iota = 1$ if $N \in 1 +2 \Z$
and the action of the symmetric group $\frak{S}_N$
on $\sum_i \C \hat{i}$ is given by 
$w \hat{i} = \widehat{w(i)}$.
Note that $(\lambda\, |\, \mu)^{\sim} \in \Z$ for every
$\lambda, \mu \in \bigoplus_i \Z \Lambda_i$. 
\begin{lem}
 Let $\zeta$ be a primitive $2 N(N + L)$-th root of unity
 and $t = \zeta^N$.
 Then the matrix $S =   S^{\iota}(\zeta)$
 is both symmetric and invertible, and satisfies 
 Verlinde's formula \eqref{Verlinde}.
 Moreover, we have:
 \begin{equation}
  \label{Szl0} 
  S^{\iota}(\zeta)_{\lambda 0} = 
  \iota^{| \lambda |} D(\lambda)_t, 
 \end{equation}
 \begin{equation}		
 \label{Szqr}
  S^{\iota}(\zeta)_{\Lambda_q \Lambda_r} =
  \iota^{q + r} \sum_s
  (\zeta t)^{-2qr} t^{2s (q + r -s + 1)}
  D ( \Lambda_{q+r-s} + \Lambda_{s} )_t.
 \end{equation}
for each $\lambda \in \V_{NL}$ and $0 < r \leq q < N$, 
where the summation in \eqref{Szqr} is taken over 
$\max \{0, q + r -N \} \leq s \leq r$.
\end{lem}
\begin{pf}
 Suppose $\zeta = \zeta_0$. Then the formula \eqref{Szqr}
 follows from \eqref{SlmC} for 
 $\cal{C} (\frak{sl}_N, N + L)$, \eqref{NlLmm}
 and \eqref{dimqU}. In this case, the other assertions also follow from
 the results for  $\cal{C} (\frak{sl}_N, N + L)$.
 For other $\zeta$, the assertions follow from Galois theory
 for $\Q (\zeta_0)/ \Q$. 		
\end{pf}
\begin{lem}
Let $N_{\lambda \mu}^{\nu}$ be the fusion rules of 
$SU(N)_L$-WZW models and
let $S$ and $S^{\prime}$ be symmetric matrices 
whose entries are indexed by $\V = \V_{NL}$.
If these satisfy Verlinde's formula \eqref{Verlinde},
$ S_{\lambda 0}$ $=$ $S^{\prime}_{\lambda 0} \ne 0$
$(\lambda \in \V)$ and 
\begin{equation}		
\label{Sqr=Sqr}
\quad
 S_{\Lambda_q \Lambda_r} =
 S^{\prime}_{\Lambda_q \Lambda_r}
 \quad (0 < r \leq q < N), 
\end{equation}
then we have $S = S^{\prime}$.
\end{lem}

\begin{pf}
We recall that there exists an algebra
surjection from 
$\Z [x_1, \ldots, x_N]^{\frak{S}_N}$ onto 
$\EuScript{F}$ (cf. Theorem \ref{SSrep} (2)), 
which sends the Schur function 
$s_{(\lambda_1, \ldots, \lambda_N)}$
(see e.g. \cite{Macdonald}) to
$[ L_{\lambda} ]$ for each $\lambda \in \V$
(see e.g. \cite{GoodmanNakanishi}).
For each 
$\xi = (\xi_1, \ldots, \xi_m)$ such that  
$N > \xi_1 \geq \ldots \geq \xi_m > 0$,
we define $E_{\xi} \in \EuScript{F}$
to be the image of the elementary symmetric function
$e_{\xi}$ via this map, that is
$E_{\xi} = [\Omega^{\xi_1}] \cdots [\Omega^{\xi_m}]$.
Since $\{ e_{\xi} \}$ is a basis of 
$\Z [x_1, \ldots, x_N]^{\frak{S}_N}$, 
$ \{ E_{\xi} \}$ spans $\EuScript{F}$.
We define the symmetric bilinear forms 
$\cal{S}$ and $\cal{S}^{\prime}$ on
$\EuScript{F}$ by setting
\begin{equation}
 \cal{S}\!
 \left( [L_{\lambda}],\, [L_{\mu}] \right) 
 = 
 S_{\lambda \mu},
\quad
 \cal{S}^{\prime}\!
 \left( [L_{\lambda}],\, [L_{\mu}] \right) 
 = 
 S^{\prime}_{\lambda \mu}.
\end{equation}
Then, Verlinde's formula for $S$ is rewritten as		
\begin{equation}
\label{Verlinde2}
 \cal{S}\!
 \left(
 ab,\, \frac{[L_{\nu}]}{S_{\nu 0}} 
 \right) =
 \cal{S}\! \left( 
 a,\, \frac{[L_{\nu}]}{S_{\nu 0}} 
 \right) 
 \cal{S}\! \left( 
 b,\, \frac{[L_{\nu}]}{S_{\nu 0}} 
 \right), 
\end{equation}
where $\nu \in \V$
and $a,b \in \EuScript{F}$. 
By \eqref{Sqr=Sqr} and this formula, we obtain
\begin{equation}
 \cal{S}
 (E_{\xi},\, [\Omega^r]) =
 \cal{S}^{\prime}
 (E_{\xi},\, [\Omega^r]) 
\end{equation}
for each $\xi$ and  $0 \leq r < N$, or equivalently, 
we obtain
\begin{equation}
 \cal{S}
 ([L_{\lambda}],\, [\Omega^r]) =
 \cal{S}^{\prime}
 ([L_{\lambda}],\, [\Omega^r])
\end{equation}
for each $\lambda \in \V$ and $0 \leq r < N$.
Repeating similar consideration, we conclude that 
$S_{\lambda \mu}$ $=$
$S^{\prime}_{\lambda \mu}$
holds for every $\lambda, \mu  \in \V$. 
\end{pf}
For $\frak{S} = \SSi$, we denote the semisimple 
ribbon category $\bold{Com}^f_{\frak{S}}$ 
(resp. unitary ribbon category $\bold{Com}^{fu}_{\frak{S}}$)
by
$\cal{C}_{\frak{S}}(A_{N-1}, t)_{\epsilon, \zeta}^{\iota}$ 
(resp.
$\cal{C}^u_{\frak{S}}(A_{N-1}, t)_{\epsilon, \zeta}$). 
\begin{thm}
\label{ComSSisMTC} 
Let $N \geq 2$ and $L \geq 1$ be integers, 
$\iota = \pm 1$ if $N \in 2 \Z$, 
$\iota = 1$ if $N \in 1 + 2 \Z$ and 
$\epsilon = \pm 1$. 
Let $\zeta$ be a primitive $2 N (N + L)$-th root of unity.
\\
\rom{(1)} Suppose $N$ is odd or $\epsilon = 1$ and set
$t = \zeta^N$. Then the category 
$\cal{C}_{\frak{S}}(A_{N-1}, t)_{\epsilon, \zeta}^{\iota}$ 
is a modular tensor category 
with S-matrix $S^{\iota}(\zeta)$.
\\
\rom{(2)} Suppose $N$ is even, $\epsilon = - 1$ and
$t\,:= - \zeta^N$ is a primitive $2 (N + L)$-th root of unity.
\rom{(}Note that this implies $L \in 2 \Z$ \rom{)}.
Then the category 
$\cal{C}_{\frak{S}}(A_{N-1}, t)_{-1, \zeta}^{\iota}$ 
is a modular tensor category
with S-matrix $S^{- \iota}(\zeta)$.
\end{thm}

\begin{pf}
We will prove Part (2).
By \eqref{Sl0C}, \eqref{d=D}, \eqref{d=-D}, 
\eqref{Szl0} and 
$D(\lambda)_{-t} = (-1)^{| \lambda |} D(\lambda)_t$, 
we have 
$S_{\lambda 0}^{\frak{S}}$ $=$ 
$S^{- \iota} (\zeta)_{\lambda 0}$.
Hence by the lamma above, it suffices to show that
$S_{\Lambda_q \Lambda_r}^{\frak{S}}$ $=$ 
$S^{- \iota} (\zeta)_{\Lambda_q \Lambda_r}$
for each $0 < r \leq q < N$.
As we will see in the next section, the action of
$c_{\Omega^r \Omega^q} \circ c_{\Omega^q \Omega^r}$ 
on the one-dimensional space 
$( \Omega^q \bar{\otimes} \Omega^r)
(0,\, \Lambda_{q+r-s} + \Lambda_{s})$
$=$ 
$\C \omega_q (0,\, \Lambda_q) \otimes
\omega_r (\Lambda_q,\, \Lambda_{q+r-s} + \Lambda_{s})$
is given by the scalar 
$(\zeta t)^{-2qr} t^{2s (q + r -s + 1)}$.
Hence, by \eqref{Trqf}, we obtain
\begin{multline}
\label{SLLS=Tr}
 S_{\Lambda_q \Lambda_r}^{\frak{S}} =
 \sum_s		
 \Tr \bigl( \cal{M}_{\iota} \circ c_{\Omega^r \Omega^q} 
 \circ c_{\Omega^q \Omega^r} 
 \!\!
 \bigm|_{( \Omega^q \bar{\otimes} \Omega^r)
 (0,\, \Lambda_{q+r-s} + \Lambda_{s})}
 \bigr) \\
 = (- \iota)^{p + q} \sum_s		
 (\zeta t^{\prime})^{-2qr} {t^{\prime}}^{2s (q + r -s + 1)}
 D(\Lambda_{q+r-s} + \Lambda_{s})_{t^{\prime}}, \\
\end{multline}
where $t^{\prime} = \zeta^N$
and the summation is taken over 
$\max \{0, q + r -N \} \leq s \leq r$.
Since the right-hand side of \eqref{SLLS=Tr}
equals $S^{- \iota} (\zeta)_{\Lambda_q \Lambda_r}$
by \eqref{Szqr}, this completes the proof of Part (2). 
\end{pf}
\begin{cor}
 Let $N \geq 2$ and $L \geq 1$ be integers, 
 $\epsilon = \pm 1$ and
 $t = \exp (\pm \frac{\pi i}{N + L})$.
 If $N + L \in 1 + 2\Z$, 
 $\cal{C}^u_{\frak{S}}(A_{N-1}, t)_{\epsilon, \zeta}$
 is a unitary MTC provided that $\epsilon = 1$
 or $N$ is odd, where $\zeta$ denotes an arbitrary primitive
 $2 N (N + L)$-th root of unity such that $\zeta^N = t$.
 If $N + L \in 2\Z$, 
 $\cal{C}^u_{\frak{S}}(A_{N-1}, \pm t)_{\epsilon, \zeta}$
 is a unitary MTC
 for each primitive $2 N (N + L)$-th root $\zeta$ of unity 
 such that $\zeta^N = \pm \epsilon^{N-1} t$.
\end{cor}							
\begin{rem}		
 \rom{(1)} When $N \in 1 + 2 \Z$, 
 $\frak{S}(A_{N-1},t)_{-1, \zeta}^{\iota}$
 is isomorphic to a 2-cocycle deformation of 
 $\frak{S}(A_{N-1},t)_{1, \zeta}^{\iota} $. 
 Hence 
 $\cal{C}_{\frak{S}}(A_{N-1}, t)_{1, \zeta}^{\iota}$ and 
 $\cal{C}_{\frak{S}}(A_{N-1}, t)_{-1, \zeta}^{\iota}$ 
 are equivalent. 
 
 \rom{(2)} For $\frak{g}$ $=$ $\frak{so}_N$ and $\frak{sp}_N$, 
 a category-theoretic construction of unitary MTC's 
 related to $\cal{C}(\frak{g}, \kappa)$ is
 given by Turaev and Wenzl \cite{TuraevWenzl2}. 
\end{rem}

\section{Braidings on $\Omega$}


In this section, we give some explicit
calculation of the braiding
$c_{q, r} := c_{\Omega^q \Omega^r}$ 
in order to complete the proof of 
Lemma \ref{cformula} and Theorem \ref{ComSSisMTC}.   					
Since the braiding is a natural transformation
and the multiplication of $\Omega$
gives a $\SS$-comodule map 
$m_{q, r}\!: 
\Omega^q \bar{\otimes} \Omega^r \to \Omega^{q+r}$,
$c_{q, r}$ satisfies	
\begin{equation}
\label{cOq+q'Or}
 c_{q+q^{\prime}, r} \circ
 (m_{q, q^{\prime}} \bar{\otimes} \mathrm{id}_{\Omega^r}) =
 (\mathrm{id}_{\Omega^r}
 \bar{\otimes} m_{q, q^{\prime}}) \circ
 (c_{q, r} \bar{\otimes} 
 \mathrm{id}_{\Omega^{q^{\prime}}}) \circ
 (\mathrm{id}_{\Omega^q} \bar{\otimes}
 c_{q^{\prime}, r}),
\end{equation}		
\begin{equation}
\label{cOqOr+r'}
 c_{q, r+r^{\prime}} \circ
 (\mathrm{id}_{\Omega^q} \bar{\otimes} m_{r, r^{\prime}})
 = 
 (m_{r, r^{\prime}} \bar{\otimes} 
 \mathrm{id}_{\Omega^q}) \circ
 (\mathrm{id}_{\Omega^r} 
 \bar{\otimes} c_{q, r^{\prime}}) \circ
 (c_{q, r} 
 \bar{\otimes} \mathrm{id}_{\Omega^{r^{\prime}}}).
\end{equation}		
\begin{lem}
For each $1 \leq p < N$ and $1 \leq q \leq N-p$, 
we have 
\begin{multline}
\label{cOqO1}
 c_{q, 1}  
 \bigl( 
 \omega \left( \Lambda_p\, |\, p+1,
 \ldots, p+q \right) 
 \otimes \omega \left( \Lambda_{p+q}\, |\, 1 \right) 
 \bigr) \\
 = - (- \zeta)^{-q} t^{p+1}
 \frac{[q]}{[p+1]} \,
 \omega \left( \Lambda_p\, |\,  p+1 \right) 
 \otimes 
 \omega \left( \Lambda_{p+1}\, |\,
 p+2, \ldots, p+q, 1 \right) \\
 + (- \zeta)^{-q} (- \epsilon)^q
 \frac{[p+q+1]}{[p+1]} \,
 \omega \left( \Lambda_p\, |\, 1 \right) 
 \otimes 
 \omega \left( \Lambda_{p} + \hat{1}\, |\,
 p+1, \ldots, p+q \right). 
\end{multline}
\end{lem}
\begin{pf}
 Suppose \eqref{cOqO1} is valid 
 for each $p$ and for some $q < N-p$.
 Then, by \eqref{cOq+q'Or}, we obtain
 \begin{equation}
 \label{cOq+1O1}
  c_{1+q, 1}  
  \bigl( 
  \omega \left( \Lambda_p\, |\, p+1,
  \ldots, p+q+1 \right) 
  \otimes \omega \left( \Lambda_{p+q+1}\, |\, 1 \right) 
  \bigr) \qquad\qquad\qquad
 \end{equation}
 \begin{multline*}
  = \;
  (\mathrm{id}_{\Omega_1} \bar{\otimes} m_{1, q})
  \circ
  (c_{1, 1} \bar{\otimes} \mathrm{id}_{\Omega^q})
  \Bigl[ \Bigr.
  \omega (\Lambda_p | p+1) \otimes \\
  \bigl\{ \bigr.
  - (- \zeta)^{-q} t^{p+2}
  \frac{[q]}{[p+2]} \,
  \omega \left( \Lambda_{p+1}\, |\,  p+2 \right) 
  \otimes 
  \omega \left( \Lambda_{p+2}\, |\,
  p+3, \ldots, p+q+1, 1 \right) \\
  + (- \zeta)^{-q} (- \epsilon)^q
  \frac{[p+q+2]}{[p+2]} \,
  \omega \left( \Lambda_{p+1}\, |\, 1 \right) 
  \otimes 
  \omega \left( \Lambda_{p+1} + \hat{1}\, |\,
  p+2, \ldots, p+q+1 \right)
  \bigl. \bigr\}  
  \Bigl. \Bigr] 
 \end{multline*}
 \begin{multline*}
  = \;
  \biggl\{ \biggr.
  - 
  t^{p+2}
  \frac{[q]}{[p+2]} \,
  w_{N,t}\!
   \begin{bmatrix}
    \Lambda_p             & \Lambda_{p+1} \\
    \Lambda_{p+1} & \Lambda_{p+2} 
   \end{bmatrix} 
  + 
  \frac{[p+q+2]}{[p+2]} 
  w_{N,t}\!
   \begin{bmatrix}
    \Lambda_p             & \Lambda_{p+1} \\
    \Lambda_{p+1} & \Lambda_{p+1} + \hat{1} 
   \end{bmatrix}
  \biggl. \biggr\} \\ 
  (- \zeta)^{-q} 
  \omega (\Lambda_p | p+1) \otimes 
  \omega \left( \Lambda_{p+1}\, |\,  p+2,
  \ldots, p+q+1, 1 \right) \\
  + (- \zeta)^{-q} (- \epsilon)^q
  \frac{[p+q+2]}{[p+2]} \,
  w_{N,t}\!
   \begin{bmatrix}
    \Lambda_p           & \Lambda_{p+1} \\
    \Lambda_p + \hat{1} & \Lambda_{p+1} + \hat{1} 
   \end{bmatrix} \\
  \omega \left( \Lambda_p\, |\, 1 \right) 
  \otimes 
  \omega \left( \Lambda_p + \hat{1}\, |\,
  p+1, \ldots, p+q+1 \right) .
 \end{multline*}
 Computing the right-hand side of the above equality, 
 we obtain \eqref{cOqO1} for $q+1$.
\end{pf}
Using \eqref{cOqO1} for $p = 1$, $q = N-1$
together with \eqref{cOq+q'Or}, we obtain
\begin{equation}
 c_{N, 1}  
 \bigl( 
 \omega \left( 0 |\, 1,
 \ldots, N \right) 
 \otimes \omega \left(0\, |\, 1 \right) 
 \bigr)
 = 
 - (- \zeta)^{- N} t\,
 [N]\, \omega \left( 0\, |\, 1 \right) 
 \otimes 
 \omega \left( 1\, |\, 2, \ldots, N, 1 \right). 
\end{equation} 
This shows that the constant $\vartheta$ in 
\eqref{cONO1=theta} equals
$\epsilon^{N-1} \zeta^{- N} t$  
and completes the proof of 
Lemma \ref{cformula}.
\begin{lem}
 We have the following relations\rom:
 \begin{multline}
 \label{cOqO1*}
  c_{q, 1} 
  \Bigl( 
  \omega \left( \Lambda_p\, |\, p+1,
  \ldots, p+s, 1, \ldots, q-s \right) 
  \otimes \omega \left( 
  \Lambda_{p+s} + \Lambda_{q-s}\, |\, q-s+1 \right) 
  \Bigr) \\
  \in - (- \zeta)^{-q} t^{p-q+s+1}
  \frac{[s]}{[p+1]}\, 
  \omega \left( \Lambda_p\, |\, p+1 \right) 
  \otimes 
  \omega \left( \Lambda_{p+1}\, |\,
  p+2 \ldots, p+s, 1, \ldots, q-s+1 \right) \\
  + 			
  \Omega^1 \left( \Lambda_p\, |\, 1 \right) 
  \otimes 
  \Omega^q \left( \Lambda_{p} + \hat{1},
  \Lambda_{p+s} + \Lambda_{q-s+1} \right) \\
  (1 \leq p < N, 1 \leq s \leq N-p, s < q \leq p+2s-1),
 \end{multline}
 \begin{multline}
 \label{cOqOr}
  c_{q, r} 
  \Bigl( 
  \omega \left( \Lambda_p\, |\, p+1,
  \ldots, p+q \right) 
  \otimes \omega \left( 
  \Lambda_{p+q}\, |\, 1, \ldots, r \right) 
  \Bigr) \\
  \in (- 1)^r (- \zeta)^{- qr} t^{pr+r}
  \frac{[q]! \, [p]!}{[p+r]! \, [q-r]!} \,
  \omega \left( \Lambda_p\, |\, p + 1, \ldots, p + r \right) 
  \qquad\qquad \\
  \qquad\qquad \otimes 
  \omega \left( \Lambda_{p+r}\, |\,
  p+r+1, \ldots, p+q, 1, \ldots, r \right) \\
  \qquad\qquad\qquad\qquad + 
  \sum_{\lambda \ne \Lambda_{p+r}}
  \Omega^r
  \left( \Lambda_p, \lambda \right) 
  \otimes 
  \Omega^q \left( \lambda,
  \Lambda_{p+q} + \Lambda_{r} \right) \\
  (0 \leq p \leq N-1, 0 \leq q \leq N-p, 0 \leq r \leq q),
 \end{multline}
 where $[n]! = [n] \cdots [2][1]$ and $[0]! = 1$. 
\end{lem}
\begin{pf}
 The relation \eqref{cOqO1*} follows from
 \eqref{cOq+q'Or}, \eqref{cOqO1} and 
 \begin{multline}
  c_{q-s, 1}  
  \bigl( 
  \omega \left( \Lambda_{p+s}\, |\, 1,
  \ldots, q-s \right) 
  \otimes \omega \left( \Lambda_{p+s} + \Lambda_{q-s}\, 
  |\, q-s+1 \right) 
  \bigr) \\
  = (- \zeta t)^{s - q}
  \omega \left( \Lambda_{p+s}\, |\,  1 \right) 
  \otimes 
  \omega \left( \Lambda_{p+s} + \hat{1}\, |\,
  2, \ldots, q-s+1 \right). 
 \end{multline}
 The relation \eqref{cOqOr} is easily proved  
 by induction on $r$, using 
 \eqref{cOqOr+r'} and \eqref{cOqO1*}.
\end{pf}
Using \eqref{cOqOr}, \eqref{cOqOr+r'} and 
\begin{multline}
 c_{q, r-s} \bigl(
 \omega \left( 0\, |\, 1, \ldots, q \right)
 \otimes \omega \left( \Lambda_q\, |\, q+1,
 \ldots, q+r-s \right) \bigr) \\
 = (- \zeta t )^{q s - q r}			
 \omega \left( 0\, |\, 1, \ldots, r - s \right)
 \otimes \omega \left( \Lambda_{r-s}\, |\, r-s+1,
 \ldots, q+r-s \right), 
\end{multline}
we obtain the following.

\begin{lem}
 For each $0 \leq q, r < N$ and
 $\max \{ 0, q+r-N \} \leq s \leq \min \{q,r\}$,
 we have
 \begin{multline}
  c_{q, r}
  \bigl(
  \omega \left( 0\, |\, 1, \ldots, q \right)
  \otimes \omega \left( \Lambda_q\, |\, q+1,
  \ldots, q+r-s, 1, \ldots s
  \right)
  \bigr) \\
  = \; (-1)^s (- \zeta t)^{- qr} t^{s(q+r-s+1)}
  \frac{[q]!\,[r-s]!}{[r]!\,[q-s]!} \,
  \qquad \qquad\qquad \\ \qquad\qquad
  \omega \left( 0\, |\, 1, \ldots, r \right)
  \otimes \omega \left( \Lambda_r\, |\, r+1,
  \ldots, q+r-s, 1, \ldots, s
  \right). \\
 \end{multline}
\end{lem}
As an immediate consequence of the lemma above, 
we see that $c_{r,q} \circ c_{q,r}$ acts on 
$( \Omega^q \bar{\otimes} \Omega^r)
(0,\, \Lambda_{q+r-s} + \Lambda_{s})$
as the scalar 
$(\zeta t)^{- 2 qr} t^{2 s(q+r-s+1)}$.
Thus we complete the proof of Theorem \ref{ComSSisMTC}.

\section{ABF models and $SU(2)_L$-SOS algebras}

%
%
%
In this section, we give an explicit description of 
the representation theory
of $\Ss$. 
We identify $\G = \G_{2,L}$ with the
Dynkin diagram of type $A_{L+1}$:
\begin{equation}
 \begin{matrix} 
  & & & & & & & \\ 
  \scriptstyle{0} & & \scriptstyle{1} & & & 
  \scriptstyle{L - 1} & & \scriptstyle{L} \quad\,\, \\
  \circ & \overrightarrow{ \longleftarrow \!\!\!-\!\!\!-\!\!\!-\!\!\!-} &
  \circ & \overrightarrow{ \longleftarrow \!\!\!-\!\!\!-\!\!\!-\!\!\!-} &
  \cdot \quad \cdot \quad \cdot &
  \circ & \overrightarrow{ \longleftarrow \!\!\!-\!\!\!-\!\!\!-\!\!\!-} &
  \circ \quad .
 \end{matrix} 
\end{equation}
%
%
Also, we identify $\V$ and $\G_{ij}^k$ with
$\left\{ 0, 1, \cdots, L \right\}$ and 
\begin{equation}
 \left\{ \left. (i_0, i_1, \cdots, i_k) \right| \,  
 0 \leq i_0, \cdots, i_k \leq L, \,\,
 \left| i_{\nu} - i_{\nu - 1} \right| = 1 \,
 (1 \leq \nu \leq k) \right\}
\end{equation}
respectively.
We define the set $\cal{B}$ by
\begin{equation}
 \cal{B} = \left\{ \left. \binom{k}{ij} \right| \,  
 {{i, j, k \in \V, \, |i - j | \leq k \leq i + j,}
 \atop {i + j + k \in 2 \Bbb{Z}, \, i + j + k \leq 2L}} 
 \right\}.
\end{equation}
Then, we have
\begin{equation}
 N_{ij}^k = 
 \begin{cases}
 1 & \left( \binom{k}{ij} \in \cal{B} \right) \\
 0 & \left( \text{otherwise} \right)
 .
 \end{cases}
\end{equation}
In order to simplify the formula for quantum
invariants stated in the introduction, 
we use the rational basis of 
type $\Sigma$ instead of type $\Omega$
(cf. Remark \ref{rembase}). 
The corresponding Boltzmann weight $w = w_{2,t,\epsilon}^{\Sigma}$
is given by
\begin{equation}
 w
 \begin{bmatrix}
  i & i \pm 1 \\
  i \pm 1 & i 
 \end{bmatrix}
 = - \zeta^{-1}\, \frac{\pm t^{\mp (i + 1)}}{[i + 1]},
\quad
 w 
 \begin{bmatrix}
  i & i \pm 1 \\
  i \mp 1 & i 
 \end{bmatrix}
 = \zeta^{-1} \epsilon\, \frac{[i + 1 \pm 1]}{[i + 1]},
\end{equation}
\begin{equation}
 w 
 \begin{bmatrix}
  i & i \pm 1 \\
  i \pm 1 & i \pm 2 
 \end{bmatrix}
 = \zeta^{-1} t,
\quad
 w \Bigl[ \text{otherwise} \Bigr] = 0.
\end{equation}
%
%
%
%
%
Next, we recall a realization of $L_k$ introduced in \cite{subf}.
Let $\Sigma$ be an algebra generated by the symbols
$\sigma (\p)$ $(\p \in \G^k, k \geq 0)$ 
with defining relations:
\begin{equation}
 \sumk \sigma (k) \, = \, 1,
\end{equation}
\begin{equation}
 \sigma (\p) \sigma (\q) \, = \,
 \delta_{\en (\p) \st (\q)} \, \sigma (\p \cdot \q),
\end{equation}
\begin{equation}
 \qquad
 \sigma (i, i+ 1, i) = \epsilon \sigma (i, i - 1, i)  
 \quad (0 < i < L),
\end{equation}
\begin{equation}
 \sigma (0, 1, 0) = \sigma (L, L - 1, L) = 0.  
\end{equation}
We define the grading
$\Sigma = \bigoplus_{k \geq 0} \Sigma^k$ via 
$\Sigma^k = 
\mathrm{span} \{ \sigma(\p) \left. \right| \p \in \G^k \}$.
Then each component $\Sigma^k$ becomes 
a right $\Ss$-comodule via 
\begin{equation}
 \sigma (\q) \mapsto 
 \sum_{\p \in \G^k} \sigma (\p) \otimes e \binom{\p}{\q}
 \quad (\q \in \G^k).
\end{equation}
For each $\binom{k}{ij} \in \cal{B}$, 
the element 
$\sigma_k (i,j)\!:= \epsilon^{\Len (\q)}\sigma (\q)$
does not depend on the choice of $\q \in \G^k_{ij}$,
where $\Len$ is as in \S 7, that is, 
\begin{equation}
 \Len (i, i - 1, \cdots, (i + j - k)/2, \cdots, j - 1, j)
 = 0, 
\end{equation}
\begin{equation}
 \Len \left( ( \cdots, n, n + 1, n, \cdots ) \right)
 = \Len \left( (\cdots, n, n - 1, n, \cdots) \right) + 1.
\label{sgndef}
\end{equation}
%
%
It is easy to see that
$\Sigma^k (i,j) = \Bbb{C} \sigma_k (i,j)$ 
for each $\binom{k}{ij} \in \cal{B}$
and that 
$\{ \left. \sigma_k (i,j) \right|
 \binom{k}{ij} \in \cal{B} \}$ 
is a linear basis of $\Sigma$.  
Since 
$\dim \Sigma^k (0,l) = \dim \Sigma^k (l,0) = \delta_{kl}$, 
we have 
$\Sigma^k \cong L_{k} \cong (\Sigma^k) \spcheck$
by Theorem \ref{SSrep} and \eqref{Mcheckij} .
More explicitly, we have the following.
\begin{prop}		
 The map $\Sigma^k \to (\Sigma^k)^{\lor};$
 $\sigma_k (i,j) \mapsto c \binom{k}{ij} 
 \sigma^{\lor}_k (i,j)$
 gives an identification of $\Ss$-comodules,
 where $\{ \sigma^{\lor}_k (j,i) \}$ denotes
 the dual basis of $\{ \sigma_k (i,j) \}$
 and the constant $c \binom{k}{ij}$ is given by
 \begin{equation}
 \label{ckijformula}
  c \binom{k}{ij}
  =
  (- \epsilon)^{(i-j)/2}							
  \frac{[(i + j + k)/2 + 1]! \, [(i - j + k)/2]! \, [(- i + j + k)/2]!} 
  {[i + 1] \, [(i + j - k)/2]!}.
 \end{equation}
 Under this identification, 
 the maps $d_{\Sigma^k}$ and $b_{\Sigma^k}$ 
 in \eqref{bMdef}-\eqref{dMdef} are given by
 \begin{equation}
 \label{dSigma}
  d_{\Sigma^k} (i) =
  \sum_j c \binom{k}{ji}^{-1} 
 \sigma_k (i,j) \otimes \sigma_k (j,i),
 \end{equation}
 \begin{equation}
  b_{\Sigma^k} 
  \left( \sigma_k (i,j) \otimes \sigma_k (j,l) \right)
  =  \delta_{il}\, c \binom{k}{ij}\, i
 \end{equation}
 respectively, where the summation in \eqref{dSigma} 
 is taken over all $j \in \V$ such that 	
 $\binom{k}{ji} \in \cal{B}$.
\end{prop}
\begin{pf}
 It suffices to show the first assertion.
 Since $\Sigma^k (i,j) = \C \sigma_k (i,j)$
 and $(\Sigma^k) \spcheck (i,j)$ $ =$
 $\C \sigma^{\lor}_k (i,j)$,
 there exists an isomorphism $\Sigma^k \cong (\Sigma^k) \spcheck$
 of the form 
 $\sigma_k (i,j) \mapsto c \binom{k}{ij} \sigma^{\lor}_k (i,j)$
 for some nonzero constant $c \binom{k}{ij}$
 $(\binom{k}{ij} \in \cal{B})$.
 To compute  $c \binom{k}{ij}$, we consider 
 the $\Ss$-right module structure on $\Sigma$ given by \eqref{wa}.
 Similarly to Lemma \ref{actOmega}, we obtain 
 \begin{equation}
 \label{sijeij}
  \sigma_k (i, j)\, e\! \binom{i,\, i \pm 1}{j,\, j \pm 1} =
  \zeta^{- k} {\epsilon}^{(\pm i \mp j + k)/2} 
  t^{(\mp i \pm j + k)/2}  
  \frac{[ \frac{i + j \mp k}{2} + 1]}{[i + 1]}
  \sigma_k (i \pm 1, j \pm 1).
 \end{equation}
 On the other hand, by \eqref{Santipode}, we obtain
 \begin{equation}
  \label{Sformulasu2}
  S \left(
  e\! \binom{i,\, i \pm 1}{j,\, j \pm 1} \right) =
  \frac{[j + 1]}{[i + 1]}\,
  e\! \binom{j \pm 1,\, j}{i \pm 1,\, i},
 \quad
  S \left(
  e\! \binom{i,\, i \pm 1}{j,\, j \mp 1} \right) =
  - \epsilon\, \frac{[j + 1]}{[i + 1]}\,
  e\! \binom{j \mp 1,\, j}{i \pm 1,\, i}.
 \end{equation}
 Using this together with \eqref{<va,w>} and \eqref{sijeij}${}_-$,
 we obtain
 \begin{equation}
  \label{svijei+j+}
  \sigma^{\lor}_k (i,j)\,
  e\! \binom{i, i+1}{j, j+1} =
  \zeta^{- k} {\epsilon}^{(i - j + k)/2} 
  t^{(- i + j + k)/2} 
  \frac{[\frac{i + j + k}{2} + 2]}{[i + 2]}
  \sigma^{\lor}_k (i + 1, j + 1).
 \end{equation}
 By \eqref{sijeij}${}_+$ and \eqref{svijei+j+}, we obtain
 \begin{equation}
  \frac{c \binom{k}{i + 1, j + 1}}{c \binom{k}{i, j}} 
  =
  \frac{[i + 1]\, [(i + j + k)/2 + 2]}{[i + 2]\, [(i + j - k)/2 + 1]}.
 \end{equation}
 Similarly, by computing 
 $\sigma_k (i, j)\, e\! \binom{i, i \pm 1}{j, j \mp 1}$,
 we obtain 
 \begin{equation}				
  \frac{c \binom{k}{i \pm 1, j \mp 1}}{c \binom{k}{i, j}} 
  = - \epsilon\,
  \frac{[i + 1]\, [(\pm i \mp j + k)/2 + 1]}{[i + 1 \pm 1]\, [(\mp i \pm j + k)/2]}.
 \end{equation}					
 By solving these recursion relations under some initial condition,
 we get \eqref{ckijformula}.
\end{pf}
\begin{prop}
\rom{(1)} The braiding						
$c\!: \Sigma^{m} \bar{\otimes} \Sigma^{n} \tilde{\to} 
\Sigma^{n} \bar{\otimes} \Sigma^{m}$
and its inverse
are given by 
\begin{equation}
 c^{\pm 1} \left( \sigma_m (h,i) \otimes \sigma_n (i,k) \right)
 =
 \sum_j w^{\pm}_{mn} \hijk{h}{i}{j}{k} 
 \sigma_n (h,j) \otimes \sigma_m (j,k),
\label{wmndef} 
\end{equation}
\begin{equation}
\label{wmnformula}
 w^{\pm}_{mn} \hijk{h}{i}{j}{k} 
 := \sum_{\bold{u} \in \G_{hj}^n} \sum_{\bold{v} \in \G_{jk}^m}
 \epsilon^{\Len (\bold{u}) + \Len (\bold{v}) + \Len (\q) + \Len (\r)}
 w^{\pm} \!\! 
 \left[ \bold{u} \,\, \frac[0pt]{\q}{\bold{v}} \,\, \r \right],
\end{equation}
where $\q$ and $\r$ denote arbitrary elements of
$\G^m_{hi}$ and $\G^n_{ik}$ respectively,  
and the summation in \eqref{wmndef} is  
taken over all $j \in \V$ such that 
$\binom{n}{hj}, \binom{m}{jk} \in \cal{B}$. \\
\rom{(2)} The ribbon functional of 
$\frak{S}(A_1, t)_{\epsilon, \zeta}^{\iota}$
acts on $L_i$ as the scalar $\theta_i^{-1}$
given by
\begin{equation}
 \label{thetasu2}
 \theta_i = \iota^i \zeta^{i (i + 2)}.
\end{equation}
\end{prop}
\begin{pf}
 Since the braiding is a natural transformation and 
 the map $\C \G^m \to \Sigma^m$; $\p \mapsto \sigma (\p)$
 is a $\Ss$-comodule map, we have 
 \begin{equation}
  c^{\pm 1} \left( \sigma_m (h,i) \otimes \sigma_n (i,k) \right)
  =
  \sum_j 
  \sum_{\bold{u} \in \G_{hj}^n} \sum_{\bold{v} \in \G_{jk}^m}
  \epsilon^{\Len (\q) + \Len (\r)} 
  w^{\pm} \!\! 
  \left[ \bold{u} \,\, \frac[0pt]{\q}{\bold{v}} \,\, \r \right]
  \sigma (\bold{u}) \otimes \sigma (\bold{v}).
\end{equation}
 Since 
 $\sigma (\p) = \epsilon^{\Len (\bold{\p})}  
 \sigma_m (i, j)$ for each $\p \in \G^m_{ij}$, 
 this proves Part (1).
 When $i = 1$, \eqref{thetasu2} follows from 
 the proof of Proposition \ref{SSrib}.
 For $i > 1$, \eqref{thetasu2} follows from 
 \eqref{SlmC} by induction on $i$.
\end{pf}

\section{The state sum invariants}

Let $L$ be a positive integer and let $\epsilon, \iota$
be elements of $\{ \pm 1\}$ such that $\epsilon = 1$ if 
$L$ is a odd integer.
Let $t$ be a primitive $2 (L + 2)$-th root of unity
and $\zeta$ a solution of $\zeta^2 = \epsilon\, t$.
Applying the general theory of TQFT to the MTC 
$\cal{C}_{\frak{S}}(A_1, t)_{\epsilon, \zeta}^{\iota}$,
we obtain an invariant $\tau_{\epsilon, \zeta}^{\iota}$
of oriented 3-manifolds.
%
%
To give an explicit description of $\tau_{\epsilon, \zeta}^{\iota}$, 
we prepare some terminologies on link diagrams.

Let $D$ be a generic link diagram in $\Bbb{R} \times (0,\,1)$ 
(viewed as a union of line segments $\overline{AB}$), 
which presents a framed link $L$ with 
components $K_1, \ldots, K_p$ (see e.g. \cite{Kassel}).
A point of $D$ is called {\em extremal} if 
the height function on $D$ attains its local maximum
or local minimum in this point, 
where the height function 
$\mathrm{ht}$ is the restriction of the projection 
$\Bbb{R} \times (0,\,1) \to (0,\,1)$ on $D$.
A point of $D$ is called {\em singular} if
it is either an extremal point or a crossing point.
We denote by $\sharp D$ the set of all singular points of $D$.
Let $\cal{E}$ be the set of all connected components of 
$D \setminus \sharp D$.
We say that $E \in \cal{E}$ {\em belongs} to $K_q$ ($1 \leq q \leq p$)
if $E$ is a subset of the image of $K_q$
via the projection $L \to D$.
Let c be a map from $\{K_1, \ldots, K_p \}$ to $\V$. 
We say that a map $\lambda\!: \cal{E} \to \cal{B};$
$E \mapsto \binom{\lambda_1 (E)}{\lambda_2 (E) \lambda_3 (E)}$
is a {\em state} on $D$ of {\em color} $c$ if  
$\lambda_1 (E) = c(K_q)$ for each component $K_q$
and $E \in \cal{E}$ belonging to $K_q$.
We denote by  $\bold{c}_{\lambda}$
the color of a state $\lambda$, and by 
$\cal{S} (D)$ the set of all states on $D$. 
Figure (A) shows a state $\lambda$ on a diagram of the Hopf link
with $6$ singular points, such that		
$\bold{c}_{\lambda}(K_1) =1$,
$\bold{c}_{\lambda}(K_2) =2$.
%
%
\setlength{\unitlength}{0.25mm}
\begin{picture}(400,380)(-70,-80)
 \put(120,0){\circle{7}}
 \put(120,240){\circle{7}}
 \put(180,60){\circle{7}}
 \put(180,180){\circle{7}}
 \put(240,0){\circle{7}}
 \put(240,240){\circle{7}}
 \put(90,110){\makebox(30,20)[l]{$\binom{2}{31}$}}
 \put(245,110){\makebox(30,20)[l]{$\binom{1}{23}$}}
 \put(125,30){\makebox(30,20)[l]{$\binom{1}{12}$}}
 \put(212,30){\makebox(30,20)[l]{$\binom{2}{02}$}}
 \put(125,192){\makebox(30,20)[l]{$\binom{1}{01}$}}
 \put(212,192){\makebox(30,20)[l]{$\binom{2}{11}$}}
 \put(30,174){\makebox(30,20)[l]{$\binom{1}{32}$}}
 \put(308,174){\makebox(30,20)[l]{$\binom{2}{L\, L-2}$}}
 \put(34,40){\makebox(30,20)[l]{$K_1$}}
 \put(312,40){\makebox(30,20)[l]{$K_2$}}
 \put(150,-60){\makebox(30,20)[l]{Figure (A)}}
 \put(0,120){\line(1,1){117}}   
 \put(0,120){\line(1,-1){117}}   
 \put(183,63){\line(1,1){57}}   
 \put(186,54){\line(1,-1){51}}   
 \put(174,66){\line(-1,1){54}}   
 \put(177,57){\line(-1,-1){54}}   
 \put(183,183){\line(1,1){54}}   
 \put(186,174){\line(1,-1){54}}   
 \put(174,186){\line(-1,1){51}}   
 \put(177,177){\line(-1,-1){57}}   
 \put(360,120){\line(-1,1){117}}   
 \put(360,120){\line(-1,-1){117}}   
\end{picture}\\
%
%

Next, we assign a complex number 
$\langle \lambda | A \rangle$ for each state
$\lambda \in \cal{S} (D)$ and
singular point
$A \in \sharp D$ as follows:				
When $(\lambda, A)$ is as in Figure (B) or Figure (C), then we set
\begin{equation}
 \langle \lambda | A \rangle 
 = c \binom{l}{hi} \delta_{hk} \delta_{ij}
 \quad \mathrm{or} \quad 
 \langle \lambda | A \rangle
 = c \binom{l}{ih}^{-1} \delta_{hk} \delta_{ij}
\end{equation}
respectively, where $c \binom{l}{hi}$ is as in 
\eqref{ckijformula}.

\setlength{\unitlength}{0.7mm}
\thicklines
%
%
%
%
%
\begin{picture}(200,95)
 \put(51.5,67.0){\line(-2,-3){20}} 
 \put(53.7,67.0){\line(2,-3){20}}  
 \put(51.4,65){\makebox(5,5)[l]{$\circ$}}
 \put(51,70){\makebox(5,5)[l]{$A$}}
 \put(27,19){\makebox(30,20)[l]{$\binom{l}{hi}$}}
 \put(70,19){\makebox(30,20)[l]{$\binom{l}{jk}$}}
 \put(46,0){\makebox(30,20)[l]{Figure (B)}}
 \put(123.8,39.9){\line(2,3){20}}   
 \put(121.8,39.9){\line(-2,3){20}}  
 \put(121.4,36.4){\makebox(5,5)[l]{$\circ$}}
 \put(121,31.4){\makebox(5,5)[l]{$A$}}
 \put(97,67){\makebox(30,20)[l]{$\binom{l}{hi}$}}
 \put(140,67){\makebox(30,20)[l]{$\binom{l}{jk}$}}
 \put(116,0){\makebox(30,20)[l]{Figure (C)}}
\end{picture}\\
\begin{picture}(200,95)
 \put(53.4,55.1){\line(2,3){16}}   
 \put(51.8,53){\line(-2,-3){16}} 
 \put(50,57.9){\line(-2,3){14.3}}  
 \put(55,50.5){\line(2,-3){14.3}}  
 \put(51.4,51.4){\makebox(5,5)[l]{$\circ$}}
 \put(44.9,51.4){\makebox(5,5)[l]{$A$}}
 \put(31,77){\makebox(30,20)[l]{$\binom{m}{hi}$}}
 \put(66,77){\makebox(30,20)[l]{$\binom{n}{jk}$}}
 \put(31,12){\makebox(30,20)[l]{$\binom{n}{cd}$}}
 \put(66,12){\makebox(30,20)[l]{$\binom{m}{ef}$}}
 \put(46,0){\makebox(30,20)[l]{Figure (D${}_{+}$)}}
 \put(125.2,57.8){\line(2,3){14.3}}   
 \put(119.8,50.8){\line(-2,-3){14.3}} 
 \put(121.8,54.9){\line(-2,3){15.7}}  
 \put(123.6,52.9){\line(2,-3){15.7}}  
 \put(121.4,51.4){\makebox(5,5)[l]{$\circ$}}
 \put(114.9,51.4){\makebox(5,5)[l]{$A$}}
 \put(101,77){\makebox(30,20)[l]{$\binom{m}{hi}$}}
 \put(136,77){\makebox(30,20)[l]{$\binom{n}{jk}$}}
 \put(101,12){\makebox(30,20)[l]{$\binom{n}{cd}$}}
 \put(136,12){\makebox(30,20)[l]{$\binom{m}{ef}$}}
 \put(116,0){\makebox(30,20)[l]{Figure (D${}_{-}$)}}
\end{picture}\\
When $(\lambda, A)$ is as in Figure (D${}_{\pm}$), we set
\begin{equation}
 \langle \lambda | A \rangle =
 \delta_{ij} \delta_{hc} \delta_{de} \delta_{fk}
 w_{mn}^{\pm} \hijk{h}{i}{e}{f},
\end{equation}
where $w_{mn}^{\pm}$ is as in \eqref{wmnformula}.

The following result follows from \S13
by a method quite similar to ``vertex models on link invariants''
(see e.g. \cite{Turaev} Appendix\, II), hence we omit the proof.

\begin{thm}
 Let $\tau_{\epsilon, \zeta}^{\iota}$
 be the invariant of closed oriented $3$-manifolds
 associated with the modular tensor category 
 $\cal{C}_{\frak{S}}(A_1, t)_{\epsilon, \zeta}^{\iota}$ 
 \rom{(}cf. \cite{Turaev}\rom{)}.
 Let $M$ be a $3$-manifold obtained by surgery on  
 $S^3$ along a framed link $L$ with $p$ components
 $K_1, \ldots, K_p$.	
 Let $D$ be a generic diagram in 
 $\Bbb{R} \times (0,\,1)$ which presents $L$.
 Then $\tau_{\epsilon, \zeta}^{\iota}$
 is given by
 \begin{equation}
 \label{tauMformula}
  \tau_{\epsilon, \zeta}^{\iota} (M) = 
  \Delta^{\sigma (L)} 
  \EuScript{D}^{- \sigma (L)- p -1} 
  \sum_{\lambda \in \cal{S} (D)} \,
  \prod_{q =1}^{p} \iota^{\bold{c}_{\lambda}(K_q)}\,
  [\bold{c}_{\lambda}(K_q) + 1] \,
  \prod_{A \in \sharp D} \langle \lambda | A \rangle.
 \end{equation}
 Here $\Delta$ denotes a fixed square root of
 $\sum_{i \in \V} [i + 1]^2$,  
 $\EuScript{D} = \sum_{i \in \V} \iota^i \zeta^{- i (i + 2)} [i + 1]^2$
 and $\sigma (L)$ denotes the signature of
 the linking matrix of $L$.
\end{thm}

\end{document}